\documentclass[11pt,reqno]{amsart}
\usepackage{amsmath, amsfonts, amssymb, amscd, enumerate, enumitem, color, ulem}
\usepackage[makeroom]{cancel}

\theoremstyle{plain}
\newtheorem{theo}{Theorem}[section]
\newtheorem{lem}[theo]{Lemma}

\theoremstyle{definition}
\newtheorem{rem}[theo]{Remark}

\newtheorem{definition}[theo]{Definition}
\newenvironment{pf}{\noindent{\it Proof. }}{$\square$\par\medskip}

\theoremstyle{plain}

\theoremstyle{definition}

\renewcommand{\=}{:=}

\newcommand{\beq}{\begin{equation}}
\newcommand{\eeq}{\end{equation}}
\renewcommand{\a}{\alpha}
\renewcommand{\b}{\beta}

\newcommand{\f}{\varphi}
\newcommand{\g}{\gamma}
\newcommand{\h}{\eta}

\renewcommand{\r}{\rho}

\renewcommand{\t}{\tau}

\newcommand{\z}{\zeta}

\newcommand{\G}{\Gamma}

\newcommand{\bB}{\mathbb{B}}
\newcommand{\bC}{\mathbb{C}}

\newcommand{\bR}{\mathbb{R}}
\newcommand{\bZ}{\mathbb{Z}}

\newcommand{\bJ}{\mathbb{J}}
\newcommand{\bH}{\mathbb{H}}

\newcommand{\bO}{\mathbb{O}}

\renewcommand{\gg}{\mathfrak{g}}

\newcommand{\gk}{\mathfrak{k}}
\newcommand{\gl}{\mathfrak{l}}
\newcommand{\gm}{\mathfrak{m}}
\newcommand{\gn}{\mathfrak{n}}

\newcommand{\gp}{\mathfrak{p}}

\newcommand{\gv}{\mathfrak v}

\newcommand{\gX}{\mathfrak{X}}

\newcommand\GL{\mathrm{GL}}
\newcommand\SL{\mathrm{SL}}
\newcommand\SO{\mathrm{SO}}
\newcommand\SU{\mathrm{SU}}
\newcommand\U{\mathrm{U}}
\newcommand\Spin{\mathrm{Spin}}
\newcommand\Sp{\mathrm{Sp}}

\newcommand{\cB}{\mathcal{B}}
\newcommand{\cC}{\mathcal{C}}
\newcommand{\cD}{\mathcal{D}}

\newcommand{\cF}{\mathcal{F}}

\newcommand{\cH}{\mathcal{H}}

\newcommand{\cM}{\mathcal{M}}

\newcommand{\cO}{\mathcal{O}}

\newcommand{\cQ}{\mathcal{Q}}

\newcommand{\cS}{\mathcal{S}}
\newcommand{\cT}{\mathcal{T}}
\newcommand{\cU}{\mathcal{U}}
\newcommand{\cV}{\mathcal{V}}
\newcommand{\cW}{\mathcal{W}}
\newcommand{\cX}{\mathcal{X}}
\newcommand{\cY}{\mathcal{Y}}
\newcommand{\cZ}{\mathcal{Z}}

\newcommand{\Jst}{J_o}

\newcommand{\p}{\partial}

\renewcommand{\square}{\kern1pt\vbox
{\hrule height 0.6pt\hbox{\vrule width 0.6pt\hskip 3pt
\vbox{\vskip 6pt}\hskip 3pt\vrule width 0.6pt}\hrule height0.6pt}\kern1pt}

\DeclareMathOperator\Aut{Aut\;}

\DeclareMathOperator\Ad{Ad}
\DeclareMathOperator\ad{ad}

\DeclareMathOperator\Id{Id}

\renewcommand\Im{\operatorname{Im}}

\newcommand{\wt}{\widetilde}
\newcommand{\wh}{\widehat}

\newcommand{\be}{\begin{equation}}
\newcommand{\ee}{\end{equation}}

\def\<#1,#2>{\langle\,#1,\,#2\,\rangle}
\newcommand{\arr}{\begin{array}{rlll}}
\newcommand{\ea}{\end{array}}
\newcommand{\bea}{\begin{eqnarray}}
\newcommand{\eea}{\end{eqnarray}}
\newcommand{\bean}{\begin{eqnarray*}}
\newcommand{\eean}{\end{eqnarray*}}

\def\sideremark#1{\ifvmode\leavevmode\fi\vadjust{
\vbox to0pt{\hbox to 0pt{\hskip\hsize\hskip1em
\vbox{\hsize3cm\tiny\raggedright\pretolerance10000
\noindent #1\hfill}\hss}\vbox to8pt{\vfil}\vss}}}

\newcounter{ssig}
\newcounter{ttig}
\newcommand{\gxx}{g(x)}
\newcommand{\gxy}{g(y)}
\newcommand{\txx}{t(x)}
\newcommand{\txy}{t(y)}

\newcommand{\Zxx}{Z(x)}

\title[Monge-Amp\`ere exhaustions 
]
{Monge-Amp\`ere exhaustions 
\\of almost homogeneous manifolds}\author[M. Kalka, G. Patrizio and A.  Spiro]{Morris Kalka, Giorgio Patrizio
and Andrea Spiro}


\thanks{{\it Acknowledgments}. This research was partially supported by the Project MIUR ``Real and Complex Manifolds: Geometry, Topology and  Harmonic Analysis'' and by GNSAGA of INdAM}

\begin{document}

\dedicatory{Dedicated to Ngaiming Mok for his sixtieth birthday}
\begin{abstract}  We  consider  three fundamental classes of compact almost homogenous  manifolds  and  show that  the complements  of  singular complex orbits in such manifolds  are endowed with  plurisubharmonic exhaustions satisfying  complex homogeneous Monge-Amp\`ere equations. This  extends to  a  new  family of mixed type examples  various classical  results on parabolic spaces and complexifications of symmetric spaces.  Rigidity results on complex spaces modeled on such new examples are given.  \end{abstract}

\maketitle

\setcounter{section}{0}
\section{Introduction}

Plurisubharmonic exhaustions satisfying the complex homogeneous Monge-Amp\`ere equation on a complex space appear naturally in many contexts.  Probably the first time they have been extensively considered was in Value Distribution Theory on affine algebraic varieties or, more generally, on parabolic spaces (\cite{GK, St0}).  When the exhaustions satisfy the complex homogeneous Monge-Amp\`ere equation with the least possible degeneracy --
the {\it strictly parabolic} case in the terminology of Stoll --  a natural foliation  is   associated to the exhaustion, namely the collection of complex curves  that are tangent to  the annihilators of the Levi form. The  very nice behavior of  these exhaustions in this case and the analogy with
the case of Riemann surfaces, suggested Stoll that these might be instrumental for the characterization of special complex manifolds such as $\bC^{n}$,  the unit ball $\bB^{n}$,  bounded complete circular domains  or affine cones. It is  well known that this is indeed the case (see for instance \cite{St1, Bu, St2, Pt1,W, PS2}). In all these instances the minimal set of the exhaustion is always a point or, after blowing up to resolve singularities, a compact projective manifold.  The exhaustion has  always a logarithmic type of singularity along  such  minimal set.  In these examples,    there is a sharp difference of behaviors, depending on whether or not the plurisubharmonic exhaustions,  satisfying the complex homogeneous Monge-Amp\`ere equation,  is bounded from above. If the  exhaustion is unbounded, the associated foliation is necessarily holomorphic, the holomorphic type of the manifold is fixed
and the only   allowed   deformations   are  rescalings  of the exhaustion in the direction of the leaves of the foliation (\cite{Bu, PS2}). On the other hand, when the exhaustion is bounded above, there is a very rich class of non biholomorphically inequivalent examples which are suitable deformation of the unit ball $\bB^{n}\subset \bC^{n}$  -- in fact an infinite dimensional class, see \cite{Le, Le2, BD, PS}   --  and, up to rescaling, the  plurisubharmonic exhaustions  that satisfy the complex homogeneous Monge-Amp\`ere equation are
 pluricomplex Green functions. \par
In a rather different development, plurisubharmonic exhaustions satisfying the complex homogeneous Monge-Amp\`ere equation occur as  a natural byproduct of   the construction of intrinsic complexifications of real analytic compact Riemannian manifolds, the so called {\it Grauert Tubes}.  They are tubular neighborhoods of  each such manifold in  its  tangent bundle -- possibly the entire tangent bundle --
 equipped with  an {\it adapted complex structure}   with the property that   the differentials of     parameterizations of geodesics of the manifold   are holomorphic embeddings   into the Grauert tube of   strip shaped neighborhoods  of the real line in ${\bC}$ (see   e.g.  \cite{GS, LS, PW}; see also \cite{DK} for generalizations to the Finsler case). Here the real analytic Riemannian manifold sits in the   Grauert tube as a top dimensional totally real submanifold, and it is the minimal set of the norm  function  on the tangent bundle  determined by the Riemannian metric.   Considering the above described adapted complex structure,  the norm function turns out  to be a plurisubharmonic exhaustion satisfying the complex homogeneous Monge-Amp\`ere equation  and exhibiting a ``square root'' type of singularity  along its minimal set, that is the real analytic Riemannian manifold on which the Grauert tube has been built.   Note that, in this case, once the size of the tubular neighborhood is fixed, whether finite or -- when possible -- infinite, there is a strongly rigid behavior: the isometric type of the minimal set of the exhaustion completely determines the complex structure of the Grauert tube.
\par
In this paper we consider a large class of almost homogeneous complex manifolds  with  cohomogeneity one actions, that is complex manifolds $M$ with  a   real Lie group $G$ of biholomorphisms  acting   with real hypersurfaces  as principal orbits. Almost homogeneous complex manifolds have been extensively studied and classified (see  \cite{HO, HS, AHR} and,  for the   strictly related topic of the classification of compact   homogeneous  CR manifolds,  \cite{AS, AS1}).
Under the additional assumption that  all  principal  $G$-orbits  are strongly pseudoconvex hypersurfaces and that  the manifold has vanishing first Betti number, it is possible  to fully describe all  compact  almost $G^\bC$-homogeneous  manifolds  with  cohomogeneity one $G$-actions  of strongly pseudoconvex principal orbits  (\cite{PoS1}).   
It turns out that, up to blow ups,  there  are three types of such manifolds (see \S  \ref{threeclasses} below for details):
 \begin{itemize}
\item {\it Type 1:}  Almost homogeneous manifolds with two compact complex manifolds as exceptional orbits; they  are all $\bC P^1$ bundles over a flag manifold.
\item {\it Type 2:}  Almost homogeneous manifolds with one compact complex exceptional orbit and one totally real;  they are the compactifications of the Morimoto-Nagano manifolds,  i.e. the compactifications of  the standard complexifications of compact symmetric spaces of rank one (CROSS). 
\item {\it Type 3:} A finite list of exceptional almost homogeneous manifolds with one compact complex exceptional orbit and one  compact exceptional orbit of mixed real/complex type; the latter  is a bundle over  a flag manifold with fiber which is  either a sphere or a real projective space of specified dimension.
\end{itemize}
Each   such  manifold  $M$   has always two singular $G$-orbits and at least one of them  is complex. If $S \subset M$ is  such  a complex orbit, then on $M_o \= M \setminus S$ there exists (see Theorem \ref{theo32} for precise statement)
a $\cC^\infty$ exhaustion
$ \t: M_o \to [0, \infty)$,    such that $\{ \t = 0\}$ is exactly the other singular $G$-orbit $S_{o}$ of $M$ and whose restriction to  $ M_{o} \setminus \{ \t = 0\}=M_{o} \setminus S_{o}$  is such that: 
\begin{itemize}
\item[1)] it  is strictly plurisubharmonic (i.e. 
$2 i \partial \overline\partial  \t = d d^c  \t > 0$);
\item[2)] there exists a smooth  function $f: (0, \infty) \to \bR$  with $d f \neq 0$ such that the composition   $ u \= f \circ  \t$   is a plurisubharmonic solution to  the  Monge-Amp\`ere equation
$(\p \overline \p  u)^n = 0$.
\end{itemize}
Furthermore if $M$ is of {\it Type 1}, then  $u$ has a logarithmic singularity along  $S_{o}$ (i.e. 
$ u \= \log \t$), while when $M$ is of {\it Type 2} or of {\it Type 3},  $u$ has a ``square root'' singularity along  $S_{o}$ (i.e. $ u \= \sqrt \t$). 
\par 
For  the almost homogeneous manifolds  of {\it Type 1} or {\it Type 2}, the existence of such exhaustion is either well known or 
not surprising. What is new in this result are the examples  given by the manifolds of   {\it Type 3} and the unified approach for the construction of the exhaustion, based on  properties of group actions   and on the   detailed analysis of the complex structure of  the Morimoto-Nagano manifolds  due to Stenzel (\cite{Ste}; see also \cite{My}). 
\par 
The main motivation for presenting in  a unified way such three families of examples   comes from the need of a common approach to the problem of deformation and of rigidity for complex manifolds  with  plurisubharmonic exhaustions that are solutions  of  the  complex homogeneous Monge-Amp\`ere equation.
There are a number of elements that   play a crucial role into    the picture,  among which we name:
\begin{itemize}
\item[a)] the nature of the minimal set and of
the singularity of the exhaustion along it;  
\item[b)] the nature of the leaves of the foliation associated to the exhaustion (either parabolic or hyperbolic Riemann surfaces) which, in turn, depends on its upper boundedness or unboundedness.
\end{itemize}
In \S 4 we start   such study  of deformations, providing the set up and defining the main tools  for such investigation, in particular the appropriate deformation tensors. Here,  we suitably extend the notions introduced in \cite{PS}  and based on
 the work of Bland and Duchamp \cite{BD}, for the analysis of deformations of the so-called {\it manifolds of circular type}, an important  family of complex manifolds that is  included in the class of  examples of {\it Type 1}. In Theorem \ref{main} we give a first result which provides informations on the deformability of examples of  {\it Type 1},  {\it Type 2} and  {\it Type 3} according to the nature of the minimal set of the exhaustion. 
\par
\smallskip
\subsection{Notation}\hfill\par
 An $n$-dimensional complex manifold $M$  is considered as a pair $(M, J)$,  
 where $J$ is the $(1,1)$ tensor field that gives the complex structure.
The operator $d^c \= d^c_J$ is defined on $k$-forms by
 $d^c = J  \circ d \circ J$,  so that $d d^c  = 2 i  \partial \overline \partial $.\par
A   CR manifold of hypersurfaces type will be    indicated as triple $(N, \cD, J)$, 
given by  a real manifold $N$ of odd dimension,  a codimension one distribution $\cD \subset TN$  and 
a smooth family  $J$  of complex structures $J_x: \cD_x \to \cD_x$, $x \in N$, 
  satisfying  the   integrability conditions  $[JX, Y] +  [X,J Y] \in \cD$ and   $[JX, JY] - [X, Y] - J[JX, Y] - J[X, JY] = 0$ for any $X, Y \in \cD$.  The holomorphic distribution of   $(N, \cD, J)$  is  the subbundle $\cD^{10} \subset T^\bC N$ of  the $+i$-eigenspaces of  $\bC$-linear maps $J_x: \cD_x^\bC \to \cD_x^\bC$. We recall that a CR manifold $(N, \cD, J)$ is Levi non-degenerate if 
  and only it the underlying real distribution $\cD$ is contact. \par
A complex space $\cX$ is  actually a pair $(\cX, \cO_{\cX})$, where $\cX$ is a Hausdorff topological space  and $\pi: \cO_{\cX} \to \cX$  is the sheaf of local $\bC$-algebras,   characterizing the complex space. Any complex manifold $M$ is identified with the complex space $(M, \cO_M)$, with $\cO_M$  sheaf of germs of  local holomorphic functions of $M$.
It is well known  (\cite{Re, Ca, Na}) that if $(\cX, \cO_{\cX})$ is a complex space carrying a $\cC^\infty$-exhaustion $\t: \cX \to [0, \infty)$,  which is strictly plurisubharmonic outside a compact set, then it always admits a  {\it Remmert reduction}, i.e. a pair $((\cY, \cO_\cY), \pi)$, formed by  a Stein space $(\cY, \cO_{\cY})$ and a proper surjective holomorphic map
 $\pi: \cX \to \cY$ such that: a)  $\pi: \cX \to \cY$ has connected fibers,  b) $\pi_*(\cO_\cX)  = \cO_\cY$ and c) for any holomorphic map $f: \cX \to \cZ$ into a Stein space $(\cZ, \cO_\cZ)$,  there is a unique holomorphic map $f': \cY \to \cZ$ such that $f = f' \circ \pi$. Geometrically,     the projection map $\pi: \cX \to \cY$  collapses all positive-dimensional compact analytic sets of $\cX$. \par
\medskip
\section{Monge-Amp\`ere  spaces}
\label{section2}
\setcounter{equation}{0}
In this paper, our interest is  focused on   the following class of complex manifolds  with  plurisubharmonic exhaustions.   
\begin{definition} \label{defMonge}  Let $\wt \cX$ be a complex manifold. A {\it  Monge-Amp\`ere $\cC^\infty$ exhaustion for $\wt \cX$} is   a $\cC^\infty$ exhaustion
$ \t: \wt \cX \to [0, T)$,  possibly with  $T = \infty$,  whose restriction to  $\wt \cX \setminus \{ \t = 0\}$  verifies: 
\begin{itemize}
\item[1)] it  is strictly plurisubharmonic (i.e. 
$2 i \partial \overline\partial  \t = d d^c  \t > 0$) 
\item[2)] there exists a smooth  function $f: (0, \infty) \to \bR$  with $d f \neq 0$ such that the composition   $ u \= f \circ  \t$   is a solution to  the  Monge-Amp\`ere equation
$(\p \overline \p  u)^n = 0$ and satisfies   the non-negativity condition  $2 i \partial \overline \partial  u  \geq 0$.
\end{itemize}
 A Stein space $\cX$   is a   {\it Monge-Amp\`ere space}  if it is the space of  a Remmert reduction $\pi: \wt \cX \to \cX$
 from a complex manifold $\wt \cX$ with  a   Monge-Amp\`ere $\cC^\infty$ exhaustion $ \t: \wt \cX \to [0, T)$. If $\cX$  is a smooth complex manifold, we call it {\it Monge-Amp\`ere manifold}. \par
 If $\cX$ is a Monge-Amp\`ere space, the  continuous function 
 $$\t': \cX  \to [0, T) \ ,\qquad \t'(x) \=   \t(y) \  \text{for some}\ y \in \pi^{-1}(x)$$
 is called {\it Monge-Amp\`ere $\cC^0$ exhaustion} of $\cX$.  The level set $\t'{}^{-1}(0) = \pi( \t^{-1}(0))$ is   the {\it soul  of $\cX$ determined by $\t$}.   By construction,  the exhaustion $\t'$ is surely  of class $\cC^\infty$ 
on the complementary set  of  the  soul. 
\end{definition}
A modeling example for the class of  Monge-Amp\`ere spaces  is   the complex Euclidean space $\bC^n$, equipped with  the 
standard  exhaustion 
\beq \label{standardex} \t_o: \bC^n \to [0, + \infty)\ ,\qquad \t_o(z) = \|z\|^2\ . \eeq
 Indeed, $\bC^n$ is the Remmert reduction of the blow up  $\pi: \wt \bC^n \to \bC^n$ of $\bC^n$ at the origin and 
the  unique smooth function $ \t: \wt \bC \to [0, + \infty)$, which extends the function  $ \t_o|_{\bC^n \setminus \{0\}} \to (0, + \infty)$  at all points of $ \pi^{-1}(0) \simeq \bC P^{n-1}$, is a  Monge-Amp\`ere $\cC^\infty$ exhaustion for $\wt \bC^n$. The function $ u = f \circ  \t$,  which  satisfies (2) in this case,  is  $ u(x) = \log( \|x\|^2)$. The soul   is  the singleton  $\{0\}$.  \par
\smallskip
 Other important examples of Monge-Amp\`ere manifolds  are    the domains of circular type (\cite{Pt, Pt1, PS}), a class which naturally includes all circular domains and all strictly convex domains of $\bC^n$. As for $\bC^n$, 
 each manifold of circular type   is the  Remmert reduction  of   its  blow-up at a fixed point $x_o$, called {\it center},  and it  has a Monge-Amp\`ere $\cC^0$-exhaustion,  whose corresponding soul  consists   only of  the center $x_o$. \par
\smallskip
Examples of Monge-Amp\`ere manifolds with  souls  containing more than  one point  are  given by the so-called {\it Morimoto-Nagano spaces}.  They are the  complex manifolds $(N, J)$, in which $N = T(G/K)$ is the    tangent bundle  of a CROSS $G/K$ and $J$  is  the  $G^\bC$-invariant complex structure, determined by the natural identification of $T(G/K) \simeq G^\bC/K^\bC$. 
A Morimoto-Nagano space $(N, J)$   is  equipped with a Monge-Amp\`ere exhaustion $\t: T(G/K) \to [0, + \infty)$ for which $\t^{-1}(0)$ coincides with 
  the zero section of $N = T(G/K)$ and is therefore   a totally real manifold of maximal dimension.  In these examples, the  exhaustion is actually  $\cC^\infty$ at all points  and  the  manifold  is  the  Remmert reduction of itself.  \par
  \medskip
\section{Monge-Amp\`ere exhaustions of   almost homogeneous manifolds}
\setcounter{equation}{0}
\subsection{Almost homogeneous manifolds  with   cohomogeneity one actions}\hfill\par
Let $M$ be an $n$-dimensional complex manifold  and   $G$  a   real Lie group of biholomorphisms of $M$, which  
acts  on $M$  {\it of cohomogeneity one}, that is  with principal orbits   that  are  
 real hypersurfaces of $M$. Since  any  principal orbit $G{\cdot} x$ is a real hypersurface of  a complex manifold, 
it is naturally equipped with a $G$-invariant induced  CR structure $(\cD, J)$ of hypersurface type. 
 If  all  principal  $G$-orbits   are strongly pseudoconvex, we say that  the cohomogeneity one action  is   {\it of strongly pseudoconvex type}. \par
 \smallskip
A   notion that is strictly related with the cohomogeneity one actions is the  following.  Let   $G^\bC $ be a {\it complex} Lie group of 
biholomorphisms of    $M$.  The complex manifold  is called {\it almost homogeneous  for  to the $G^\bC$-action}   (shortly,   {\it almost $G^\bC$-homogeneous}) if there is  a  $G^\bC$-orbit 
which is  open and dense in  $M$.\par
\smallskip
 It is clear that any homogeneous complex manifold $M = G^\bC/H$ is  almost homogeneous, but there are many examples of almost homogeneous manifolds that are not homogenous.  
Many  such  examples  are indeed offered  by  cohomogeneity one actions.  Assume that $G \subset G^\bC$ is a  compact real form of a reductive complex Lie group  $G^\bC$ and  that there is  a cohomogeneity one $G$-action  on a compact complex manifold $M$. If $x \in M$ is a regular  point for  the $G$-action, the  $G^\bC$-orbit $G^\bC{\cdot} x$ is   a complex submanifold of $M$   containing the real hypersurface $G{\cdot}x$.  It is therefore  open in $M$.  As a consequence of   standard facts    on the orbit space of cohomogeneity one actions  (see e.g. \cite{Br, PoS0, PoS1}),  one can   see   that such open orbit $G^\bC {\cdot} x$ is  dense but,  in general,   not equal to  $M$. More precisely,  $G^\bC{\cdot}x = M$ if and only if the real Lie group  $G$ has  no complex  singular orbits in $M$.  \par
\smallskip
\subsection{Three important classes  of  almost homogeneous manifolds}\hfill\par
\label{threeclasses}
We now focus  on  a special class of  compact almost homogeneous  $G^\bC$-manifolds with   a compact real form  $G \subset G^\bC$  acting of     
 cohomogeneity one.\par
\smallskip
Let $M$ be a compact complex manifold with a cohomogeneity one holomorphic $G$-action  of 
strongly pseudoconvex type. Each principal $G$-orbit  $N = G{\cdot}x \subset M$  has the following two important properties: 
\begin{itemize}
\item[--] it is a compact homogeneous $G$-manifold,  identifiable with a coset space $N = G/H$; 
\item[--] it  has an induced $G$-invariant  strongly pseudoconvex CR structure $(\cD, J)$. 
\end{itemize}
 The classification (up to coverings) of  compact   homogeneous  CR manifolds  $(G/H, \cD, J)$  with these two properties
 has been determined in  \cite{AS, AS1}. From  this classification and other  important properties of  
  almost homogeneous manifolds, proved in \cite{HO, HS, AHR},   in principle one can get a complete  description of all 
  compact  almost $G^\bC$-homogeneous  manifolds $M$ with  cohomogeneity one $G$-actions  of strongly pseudoconvex type.  
 Such  description is given explicitly    in \cite{PoS1} under the   assumption that the first Betti number is $b_1(M) = 0$. 
As it is pointed out in  \cite{AHR},  the cases  with  $ b_1(M) = p  \neq 0$ are  fibered bundle  over  $p$-dimensional complex tori, with a fibre $M'$ which is an almost homogeneous manifold    with $b_1(M') = 0$.   \par 
  \smallskip
According to  \cite{PoS1}, any compact  almost homogeneous manifold  $M$ with  strongly pseudoconvex, cohomogeneity one $G$-actions and   with $b_1(M) = 0$
belongs to one of   the following disjoint  three  classes.  Here,  the complex structure   $J$  of $M$ is  the natural  $G$-invariant complex structure.
 \par
\smallskip
\subsubsection{Almost homogeneous manifolds with two ends}
\label{3.2.1}
Consider the  $\bC P^1$-bundles   of the form
\beq \label{incan} \pi: M = G^\bC \times_{P, \r} \bC P^1 \longrightarrow  G^\bC/P\ ,\eeq
where: 
\begin{itemize}
\item[a)] $G^\bC/P$ is a flag manifold (i.e. a homogeneous quotient of a complex semisimple Lie group $G^\bC$ by a parabolic subgroup $P \subset G^\bC$), equipped with 
a  compact real form $G \subset G^\bC$ and a fixed choice of a $G$-invariant  K\"ahler metric $g$ on $G^\bC/P$; 
\item[b)] $\rho: P \to \Aut(\bC P^1)$ is a biholomorphic  action on $\bC P^1$ of the  isotropy group $P$,  such that 
  $\r|_{G\cap P}: G \cap P\to   \Aut(\bC P^1)$  coincides with the 
standard  cohomogeneity one action of $T^1 = \rho(G \cap P) $ on $\bC P^1$. 
\end{itemize}
For each such $\bC P^1$-bundle, the compact real form $G \subset G^\bC$ acts transitively on the flag manifold $F = G^\bC/P$ and of cohomogeneity one on $M$. There are  exactly two singular 
$G$-orbits, say $S$, $S'$,  both of them  complex and $G$-equivalent to the flag manifold $F = G^\bC/P = G/G\cap P$. Their  intersections  with a fibre of $\pi: M = G^\bC \times_{P, \r} \bC P^1 \longrightarrow  G^\bC/P$ are the  two singular orbits of the standard action of $T^1$ on $\bC P^1$. The action of $G^\bC$   has three disjoint orbits: 
$M_{\text{reg}} \= M \setminus (S \cup S')$ (which is  open  and dense), $S$ and $S'$.  \par
\smallskip
A  manifold of this kind is usually called {\it with two ends}, since  any  singular  $G$-orbit that is complex is referred as an  {\it end} of the manifold. 
It is known that any other almost homogeneous manifold, which  satisfies the above conditions  and for which there are  two complex singular $G$-orbits,   admits a   blow up, which is $G$-equivalent to  $\bC P^1$-bundles described  above (\cite{PoS1}, Thm. 2.4). We shortly call such manifolds  {\it almost homogeneous manifolds with two ends} and those as in \eqref{incan}   {\it  in canonical form}. \par
\smallskip
The simplest example in this class   is the blow up $\wt{\bC P^n}$ of $\bC P^n$ at a point $[x_o] $. Indeed,  it is  an almost homogeneous manifold 
 with two ends with  
 \begin{align*}
 & G = \SU_n\ ,\qquad G^\bC = \SL_2(\bC)\ ,\\
 & G^\bC/P \= \SL_n(\bC)/P \simeq \bC P^{n-1}\ ,\qquad \text{where}\  P \= \{A \in \SL_n(\bC)\ : A {\cdot} x_o  = x_o\}
 \end{align*}
  The singular $G$-orbits  are  both biholomorphic to $\bC P^{n-1}$, one given by the exceptional divisor at $[x_o]$, the other by    the hyperplane $\pi_o = \{[x]:x \in (x_o)^\perp\}$.  
\par
\smallskip
\subsubsection{Compactifications of Morimoto-Nagano spaces}  
\label{3.2.2}  It is the class   of compact complex manifolds, given by the infinite sequences of manifolds
$$\bC P^n\ ,\ \ \ \cQ^n = \{[z] \in \bC P^{n+1}\ :\ {}^tz z = 0\}\ ,\ \ \bC P^n \times \bC P^n\ ,\ \ \operatorname{Gr}_{2, 2n}(\bC)\ ,$$
together with   the  {\it Cayley projective plane} 
$$EIII = E_6/\SO_2 {\cdot} \Spin_{10}\ .$$
Each  of these manifold is a $G$-invariant {\it compactification of a Morimoto-Nagano space} for an appropriate compact simple Lie group $G$. More precisely, 
\begin{itemize}
\item[a)]  $\bC P^n$ is the  $\SO_n$-invariant complex compactification of    $T \bR P^n$, 
\item[b)]  $\cQ^n$ is the $\SO_n$-invariant complex compactification of  $T S^n$, 
\item[c)] $\bC P^n \times \bC P^n$ is the $\SO_n$-invariant complex compactification of $T \bC P^n$, 
\item[d)] $\operatorname{Gr}_{2, 2n}(\bC)$  is  the  $\Sp_n$-invariant compactifications of $T \bH P^n$, 
\item[e)]  $EIII$  is the   $\operatorname{F}_4$-invariant complex compactification of $T \bO P^2$.
\end{itemize}
In all these cases,  $G$ acts of cohomogeneity one of  strongly pseudoconvex type. There are   two singular $G$-orbits,  one  complex, the other 
totally real. The former is the complex manifold which is  complementary to  the tangent bundle $TS$ of the CROSS $S = G/K$,  the latter is   the zero section of  $TS$ and is therefore identifiable with $S = G/K$.
Since  only one singular $G$-orbit  is complex, these   manifolds are said to be  {\it  with one end}. \par
\smallskip
 We remark  that,  by the  results of  Morimoto and Nagano (\cite{MN}),  $\bC^n$,   the    unit ball $\bB^n \subset \bC^n$ and the Morimoto-Nagano spaces   
 are  the only   Stein manifolds on which there is  a biholomorphic cohomogeneity one action of strongly pseudoconvex type  for a  compact Lie group $G$. This is one of  the main reasons of interest for this class. 
\par
\smallskip
\subsubsection{Almost homogeneous manifolds with one end of mixed type} 
\label{3.2.3}
This class  consists of the almost homogeneous manifolds constructed as follows.  Let $G$ be a compact Lie group and $\wh M$ a  homogeneous $G$-bundle of the form
\beq \label{incan1} \pi: \wh M = G \times_{G_Q, \r} F \longrightarrow G/G_Q\ ,\eeq
where the basis $G/G_Q$,  the  fibre $F$ and the
representation $\rho: G_Q \to \Aut(F)$ form one of the triples listed in  Table 1.  There, the map $\rho$ is indicated  only by   
the    group $\rho(G_Q)$, which in all cases has to be considered as a group of  projective transformations   of a projective space $\bC P^s$ or $\bC P^{s+1}$,   depending   on whether    $F = \bC P^s$ or  $F = \cQ^s \subset \bC P^{s+1}$. \par
The almost homogeneous manifolds of this class  are those having the form \eqref{incan1} together with all other complex $G$-manifolds with exactly  one complex singular $G$-orbit and admitting  a manifold  $\wh M$ as a blow-up along  such $G$-orbit.  By a direct inspection of the Levi forms of the regular $G$-orbits (they can be determined from  the  explicit descriptions   in \cite{AS2}), one can check that all of  them   are  of strongly pseudoconvex type. Each manifold of this third class has two singular   $G$-orbits $S, S'$, the first  complex,   the second neither complex nor totally real. However, the intersection of $S'$  with each fibre $F_x= \pi^{-1}(x)$ (which is an almost homogeneous space of the second class)  is a totally real submanifold $G$-equivalent to  $\bR P^n$ or $S^n$.   Due to this, 
the   manifolds of this third class are called    {\it  with one end and of mixed type} and those  as in  \eqref{incan1}    {\it in canonical form}. \par
\smallskip 
 \centerline{
{\renewcommand{\arraystretch}{1.2}
\begin{tabular}{|c|c|c|c|c|}
\hline
\ &  $G/G_Q$ &  $F$ & $\rho(Q)$ \\
\hline
\hline
$I_1$ & $\SU_n/\operatorname{S}(\U_2 \times \U_{n-2})$ & $ \bC P^2$ & $\SO_3$  \\
 \hline
$I_2$ & $\SU_n/\operatorname{S}(\U_2 \times \U_{n-2})$ & $\cQ^2$ & $\SO_3$   \\
 \hline
 $II$ & ${\renewcommand{\arraystretch}{1}\begin{array}{r} \left(\SU_p/\operatorname{S}(\U_2 \times \U_{p-2})\right) \times \ \left(\SU_q/\operatorname{S}(\U_2 \times \U_{q-2})\right)\\
 p+q > 4 \end{array}}$ &   $\bC P^3$ & $\SO_4/\bZ_2$ \\
 \hline
 $III$ & ${\renewcommand{\arraystretch}{1}\begin{array}{r} \SU_n/\operatorname{S}(\U_4 \times \U_{n-4}) \qquad 
n > 4 \end{array}}$ & $\bC P^5$  & $\SO_6/\bZ_2$ \\
 \hline
 $IV_1$ & $\SO_{10}/\SO_2 \times \SO_8$ &  $\bC P^7$ & $\SO_8/\bZ_2$ \\
 \hline
 $IV_2$ &  $\SO_{10}/\SO_2 \times \SO_8$ &   $\cQ^7$ & $\SO_8$\\
 \hline
 $V_1$ & $\operatorname{E}_6/\SO_2 \times \Spin_{10}$ &  $\bC P^9$ & $\SO_{10}/\bZ_2$  \\
 \hline
$V_2$ &  $\operatorname{E}_6/\SO_2 \times \Spin_{10}$ &   $\cQ^9$ & $\SO_{10}$  \\
 \hline
 \end{tabular}
 }
 }
 \smallskip
 \centerline{\bf Table 1}
\par
\begin{rem} From Table 1,  cases $II$ and $III$ are the only ones with  no counterparts with quadrics as fibers. The reason becomes manifest   if one  recalls how   Table 1 derives from the previous  results on almost homogeneous manifolds and homogeneous CR structures. 
\par
\smallskip
By \cite{HS}, if a compact almost homogeneous $G^\bC$-manifold $M$ has a cohomogeneity one $G$-action and just one end, then it is either   a compactification of a Morimoto-Nagano space or  a fiber bundle  over a flag manifold $G/G_Q$.  In this second case,  the fiber $F$ is  either a compactification of a Morimoto-Nagano space or admits a cohomogeneity one  $G$-action with a single isolated fixed point.  Since each regular $G$-orbits in $M$ are homogeneous compact CR manifolds, 
 up to a covering,  all  of them  are further constrained to be  $G$-equivariantly equivalent to a Levi non-degenerate homogeneous CR  manifolds of the   classification in  \cite{AS}. \par
 \smallskip
  This last fact gives a lot of restrictions on the group $G$ and the fiber $F$.  In fact,  one gets that  there are very few possibilities for  the triples  $(G/G_Q, F, \rho(G_Q))$, in which $F$ must be either  $\bC P^s$ or $\cQ^s \subset \bC P^{s+1}$  with   $s = 2,3,5, 7$ or $9$,   the  group $\r(G_Q)$  either $\SO_s/\bZ$ or $\SO_s$ and   $G/G_Q$ must be  one of the five possibilities appearing  in Table I,  one per each of the  five possible cases  for the dimension $s$ of $F$. 
 An additional  restriction comes from the fact that    $G_Q$ must admit a representation $\rho: G_Q \to \Aut(F)$ with   $\rho(G_Q) = G_Q/\ker \rho$ equal to   $\SO_s/\bZ_2$ or  $\SO_s$.  For all five cases for  $s$, there exists a closed normal subgroup $N \subset G_Q$ so that $G_Q/N = \SO_s/\bZ_2$, but  only  for  three  of them there exists a  normal subgroup $N$ so that  $G_Q/N = \SO_s$.    Cases $II$ and $III$  are   those  where   there is  no  such   normal subgroup. 
   \end{rem}
  \par
\smallskip
\subsection{A new class  of Monge-Amp\`ere spaces}\hfill\par
From now on, we restrict to the complex manifolds  described in \S \ref{3.2.1}, \S \ref{3.2.2} and \S \ref{3.2.3}. For each of them we refer to  $G$ as  the    {\it  group of the cohomogeneity one action}.   \par
\smallskip
As we already mentioned,   on each   such  manifold  at least one of the two singular $G$-orbits  is complex. Let $S \subset M$ be  such  orbit and denote by $M_o \= M \setminus S$
its complementary set. By definition, it  is  a complex manifold,  on which
\begin{itemize}
\item[--] $G$ has a cohomogeneity one action of strongly pseudoconvex type; 
\item[--] $G^\bC$ acts  either transitively on $M_o$  (this occurs when the second singular $G$-orbit $S' \subset M$ is not complex) or with an open and dense orbit. 
\end{itemize}
In addition to this,  the following crucial property holds. 
\begin{theo} \label{theo32}  Let   $M$ be one of the   almost homogeneous $G^\bC$-manifolds  described   in \S\S \ref{3.2.1}, \ref{3.2.2},  \ref{3.2.3} in canonical form and  $S = G{\cdot} x\subset M$ a complex singular  $G$-orbits. Then: 
\begin{itemize}
\item[i)]   The complementary set  $M_o = M \setminus S$
   admits a $G$-invariant   Monge-Amp\`ere exhaustion $\t: M_o \to [0, \infty)$; 
\item[ii)]   This exhaustion   is the unique (up to a scaling factor) $G$-invariant Monge-Amp\`ere exhaustion $\t$ satisfying the following three conditions: 
\begin{itemize}
\item[$\a$)] it satisfies  (2) of Definition \ref{defMonge} with $f(t) = \log(t)$ in case both  singular $G$-orbits of $M = M_o \cup S$ are  complex and $f(t) = \sqrt{t}$ otherwise; 
\item[$\b$)] the level set  $\{ \t = 0 \}$ coincides with the second singular $G$-orbit $S'$ of  $M$; 
\item[$\g$)] if $S'$ is not complex, then  $u =  f \circ  \t = \sqrt{\t}$  
admits a continuous extension at each point  $x_o \in \{\t = 0\} $; if  $S'$ is complex, then for  each  $x_o \in \{\t = 0\} $, there is  a system of complex coordinates $z = (z^i)$ centered at $x_o$, in which $u =  f \circ  \t = \log \t$ has  a logarithmic singularity at $x_o$, i.e. 
$$u(z) = \log \t(z) =  \log \|z\| + O(1)\  . $$
\end{itemize}
\end{itemize}
   Thus,  any Remmert reduction $\cM_o$  of  one  such manifold  $M_o$ 
   is a Monge-Amp\`ere space. \par
 \end{theo}
We remark that this  theorem  gives a whole new  class of examples of Monge-Amp\`ere spaces $(\cX, \t)$. Indeed,  all known examples of Monge-Amp\`ere  spaces are complex $n$-dimensional manifolds with a soul  $\cS$,  which is either formed by an isolated point (it is what occurs  in a 
manifold of circular type) or a totally  real submanifold of maximal  dimension, i.e. with  $\dim_\bR \cS = \frac{1}{2} \dim_\bR \cX$ (it is the case of the Morimoto-Nagano spaces). But, by the above theorem, we see that each  manifold $M_o = M \setminus S$ determined  by a manifold $M$   in  canonical form    of the third class   has a  Remmert reduction which is {\it a Monge-Amp\`ere  space with a   soul $\cS$ that   is neither a point nor a totally real submanifold of maximal dimension}. \par
\smallskip
The proof of Theorem \ref{theo32} is based on  some  properties of  almost homogeneous spaces, which we recall in the next subsections.  We begin by  introducing some  additional notation.\par
\smallskip
\subsubsection{Notational issues}
Let $G^\bC$ be  the complexification of the semisimple Lie group $G$  and  
$\gg = Lie(G)$,    $\gg^\bC = \gg + i \gg = Lie(G^\bC)$. We denote by  $\cB$ the Cartan-Killing form of   $\gg^\bC$ and 
for any subspace $\gv \subset \gg^\bC$, we indicate by $\gv^\perp$ its $\cB$-orthogonal complement in  $\gg^\bC$. The same notation is used for the $\cB$-orthogonal complements of subspaces of $\gg$.\par
For each $X \in \gg^\bC = \gg + i \gg \subset Lie(\Aut(M_o))$, we denote by 
$\wh X \in \gX(M_o)$ the corresponding {\it infinitesimal transformation} of $M_o$, i.e. the unique complete vector field whose flow $\Phi^{\wh X}_t$, $t \in \bR$, is the 
family of diffeomorphisms
$$\Phi^{\wh X}_t: M_o \to M_o\ ,\qquad \Phi^{\wh X}_t(x) \=\exp(tX) {\cdot} x\ .$$ 
We recall that the map $\wh{\phantom{X}}: \gg^\bC \longrightarrow \gX(M_o)$  between  $\gg^\bC$ and the space of  vector fields of $M_o$ is  an injective  anti-homomorphism  of Lie algebras,  i.e. 
$[\wh X, \wh Y] = - \wh{[X,Y]}$ for all $X, Y \in \gg^\bC$.\par
\smallskip
Consider now a regular point $x_o \in M_o$ and identify $G{\cdot} x_o$ with 
the coset space $G{\cdot} x_o = G/L$ with  $x_o \simeq e {\cdot} L $.  We recall that 
the surjective linear map  
$$\imath: \gg \longrightarrow T_{x_o} (G{\cdot} x_o) =  T_{eL} G/L\ ,\qquad \imath(X) = \wh X|_{x_o}$$
induces an isomorphism between  the vector space $\gl^\perp \subset \gg$, complementary to $\gl$, and  the tangent bundle of $G{\cdot}x_o$ at $x_o$. 
In the following, we constantly use such isomorphism to identify these vector spaces. In this way,   the subspace $\cD_{x_o} \subset T_{x_o} (G{\cdot} x_o)$ and the complex structure $J_{x_o}: \cD_{x_o} \to \cD_{x_o}$ of the  CR structure $(\cD, J)$  are identified with 
\begin{itemize}
\item[--] an $\Ad_L$-invariant codimension one real subspace  $\gm \subset \gl^\perp$; 
\item[--] an $\Ad_L$-invariant  complex structure $J: \gm \to \gm$.
\end{itemize}
Choosing a unitary vector  $Z \in \gl^\perp \cap  \gm^\perp$ (it is $\Ad_L$- invariant and unique up to a sign),   we get an  $\Ad_L$-invariant  decomposition 
\beq\label{decomp}  \gg = \gl + \gl^\perp = \gl + (\gm + \bR Z)\ .\eeq
Using the fact  that $\cD$ is contact,  one can show  that 
$\gl + \bR Z = C_\gg(Z)$ (see e.g. \cite{AS2}, \S 3.1).
 \par
\smallskip
\subsubsection{Distinguished curves  in the above three classes of  almost homogeneous  manifolds} 
\label{optimalcurve}
Consider now the infinitesimal transformation $J \wh Z = \wh{i Z}$ corresponding to  $i Z \in \gg^\bC$ and let 
$$\h: \bR \longrightarrow M_o\ ,\qquad \h_t \= \exp(i t Z) {\cdot} x_o = \Phi^{J\wh{Z}}_t(x_o)\ .$$
 By Thms. 3.4 and 3.7 in \cite{Sp},    the curve $\h$ has the following crucial properties. 
\begin{itemize}
\item[(1)]  It intersects each regular $G$-orbit of $M_o = M \setminus S$; in fact, there exists a $G$-invariant K\"ahler metric $g$ on $M$, with the property that  $\h_t$ is  a reparameterization of a geodesic of $g$ orthogonal to all regular $G$-orbits. 
\item[(2)] If  $S \subset M$ is the only complex singular $G$-orbit,   then $\h$ intersects the non-complex singular $G$-orbit $S' \subset M_o$; in this case,  there is no loss of generality if  we change  the starting point of $\h$  and  assume that   $x_o = \h_0$ is in  $S'$.
\item[(3)]  If both singular $G$-orbits $S, S'$ are complex, then $\h$ intersects  neither of them;  however $\lim_{t \to +\infty} \h_t $ is either in $S$ or $S'$;   changing  $Z$ into $-Z$, we may always assume  that  $\lim_{t \to + \infty} \h_t \in S'$. 
\item[(4)] Each element   in the isotropy $G_{\h_t}$ of a regular point   $\h_t$ fixes all other points of the curve; this implies that the space  $\gl^\perp \simeq T_{\h_t} G{\cdot}\h_t$ 
 is the same for all  $G$-regular points of the curve $\h$;  there is also  a  canonical  isomorphism $\cD_{\h_t}$ between the spaces    $\cD_{\h_t}$ of the CR distributions of  the $G$-orbit $G{\cdot}\h_t$, so that they are all identifiable with a fixed vector subspace $\gm \subset \gl^\perp$,  independent on $t$. 
\item[(5)] If $M = M_o \cup S$ is an almost homogeneous manifold with one end and of mixed type, then $\h_t$ is entirely included in a single fiber of the projection $\pi: M \to G/G_Q$ over the  flag manifold $G/G_Q$ described in Table 1. 
\end{itemize}
\par
\smallskip
\subsubsection{The distinguished  curves  of  Morimoto-Nagano spaces}
\label{explicit}
Assume now that  $M = M_o \cup S$ is the compactification of a Morimoto-Nagano space, so that $M_o = TS'$ for  a CROSS $S' = G/K$. Let  $\gk = Lie(K)$ and $\gg = \gk + \gp$ the corresponding $\cB$-orthogonal decomposition of $\gg$. We recall that  $M_o = TS'$ is $G$-equivariantly identifiable with $TS' \simeq G \times_{K, \r} \gp$. with $\r(K) = \Ad_K|_{\gp}$,  and that the $G$-invariant complex structure of $TS'$ is  the pull-back $J = \f^*(\Jst)$ 
of the  complex structure $\Jst$ of $G^\bC/ K^\bC$ by means of the $G$-equivariant 
diffeomorphism 
\beq\label{3.4} \f: T S' = G \times_{K, \r} \gp \longrightarrow G^\bC/K^\bC\ ,\quad \f([(g , X)]_K) \= \exp(iX) {\cdot} gK^\bC  \ .\eeq
An explicit expression  for  $J$   has been determined by Stenzel   in \cite{Ste}  (see also \cite{My}) and it 
can be described as follow. Let 
$$\pi: G \times \gp \to G \times_{K, \r} \gp$$
be  the natural quotient map. Then, for any $(g, X) \in G \times \gp$,  the  vectors  $v \in T_{(g, X)} (G \times \gp)$ can be described as pairs   $v = (Y|_g, V)$,  with $Y \in \gg$ and $V \in \gp$.  Consequently, for  each $w \in  T_{[(g, X)]_K} (G \times_K \gp)$, there is a (non-unique) element 
$(Y^{(w)}|_g, V^{(w)}) \in T_{(g, \gp)} (G \times \gp)$ with 
$$\pi_*(Y^{(w)}|_g, V^{(w)}) = w\ .$$
By \cite{Ste},  the tensor $J_{[(g, X)]_K}$ is the unique endomorphism of  $ T_{[(g, X)]_K} G \times_{K, \r} \gp$ such that 
\beq\label{XdiY} 
\begin{split}
&   J(w) = \pi_*(- (\cT_X)^{-1}(V^{(w)})|_g + \cT_X(\ad_X(Y^{(w)\gk}))|_g\ ,\ \cT_X(Y^{(w) \gp}))\ .\\
&\text{with}\ \cT_X \= \left(\frac{\sin \ad_X}{\ad_X}\right)^{-1} \circ \cos \ad_X\ ,
\end{split}
\eeq
where  for each $E \in \gg =  \gk + \gp$, we denote by $E^\gk$,  $E^\gp$ the $\cB$-orthogonal projections  into $\gk$ and $\gp$, respectively, and  the notation $\sin \ad_X$,  $\cos \ad_X$, etc.  stand for the operators defined by power series. 
Note that for each $X$ in  $\gp$ or in $\gk$, the linear operators 
$$\frac{\sin \ad_X}{\ad_X}: \gg \longrightarrow \gg\qquad \text{and}\qquad \cos \ad_X: \gg \longrightarrow  \gg$$
 are invertible, even  and preserve $\gk$, $\gp$. So, also $\cT_X$ is invertible, even  and preserves $\gk$ and $\gp$.  
 \par
With the help of this information we may now give   an  explicit description for the curve $\h: \bR \to TS'$, described  in \S \ref{optimalcurve}.
For simplicity,    assume that  the identification $T S' \simeq G \times_{K, \r} \gp$ is done in such a way  that  the $G$-regular point $y_o = \h_1 \in M_o \setminus S'$  has  the form  
$$y_o = [(e, X_o)]_K\ \qquad \text{for some}\quad 0 \neq X_o \in \gp\ .$$
Then,   the isotropy subalgebra $\gl \= \gg_{y_o}$ of $G$ at $y_o$  is 
$\gl =  \gk \cap \ker (\ad_{X_o})$.   Moreover, since $S' = G/K$ has rank one,   we know that  $\ker(\ad_{X_o}) \cap \gp = \bR X_o$  so that the  linear map 
$$ \ad_{X_o}|_{\gp \cap (\bR X)^\perp}: \gp \cap (\bR X_o)^\perp \to \gk$$
has  trivial kernel. Consider now the vector subspaces $\gp_1$, $\gp_2$, $\gm'$  of $\gg = \gk + \gp$ defined by
$$\gp_1 \=  \gp \cap (\bR X_o)^\perp \subset \gp\ ,\qquad \gp_2 \= \ad_{X_o}(\gp_1) \subset \gk\ ,\qquad \gm' \= \gp_1 + \gp_2\ .$$ 
\begin{lem} \label{lemma33}  The subspace $\gm'$  coincides with  the $\ad_{X_o}$-invariant orthogonal complement $ \gm' = (\gl + \bR X_o)^\perp$. This also implies  that   $\ad_{X_o}(\gp_2)  = \gp_1$.
\end{lem}
\begin{pf} Since $[\gl, X_o] = 0$ we have that  
$\cB(\gl, \gp_2) = \cB(\gl, [X, \gp_1]) = - \cB([\gl, X], \gp_1) = 0$, i.e. $\gp_2 \subset \gl^\perp$. This, together with the fact that 
 $\gp_2 \subset\gk \subset (\bR X_o)^\perp$,   implies that  $ \gm'  = \gp_1 + \gp_2\subset (\gl + \bR X_o)^\perp$. 
 The equality follow  by  counting  dimensions. 
 \end{pf}
By this lemma, we have $\gl^\perp = \bR X_o + \gm'$. On the other hand,  for each  $Y_1 \in \gp_1$, $Y_2 \in \gp_2$, 
\begin{align*}
\wh{(Y_1 + Y_2)}_{[(e, X_o)]_K} &= \frac{d}{dt} \exp(t(Y_1 + Y_2)){\cdot}[(e, X_o)]_K\big|_{t = 0} = \\
& =  \pi_*\left(\frac{d}{dt} \left(\exp(t(Y_1 + Y_2)), X_o\right)\big|_{t = 0}\right)= \pi_*(Y_1|_e + Y_2|_e, 0)\ .
\end{align*}
Hence, by \eqref{XdiY}
\beq \label{Jvect}  J \wh{(Y_1 + Y_2)}_{[(e, X_o)]_K} = \pi_*( \cT_{X_o}(\ad_{X_o}(Y_2))|_e\ ,\ \cT_{X_o}(Y_1))\ .\eeq
We now observe that,  since $\cT_{X_o}$ is a series of even powers of $\ad_{X_o}$, by Lemma \ref{lemma33}, it preserves $\gp_1$ and $\gp_2$. Hence, 
 both $Y_1'\= \cT_{X_o}(\ad_X(Y_2))$ and $Y_2' \= \cT_{X_o}(Y_1)$ are in $\gp_1$. Moreover, since $\ad_{X_o}|_{\gp_2}:\gp_2 \to \gp_1$
is a linear isomorphism, there exists $Y''_2 \in \gp_2 \subset \gk$ such that $Y_1' = - [Y_2'', X_o]$. It follows that the vector \eqref{Jvect} can be written as 
\beq
\begin{split} J \wh{(Y_1 + Y_2)}_{[(e, X_o)]_K} &= \pi_*( Y'|_e\ ,\ - [Y''_2, X_o]) =
\pi_*( Y'_1|_e + Y''_2|_e\ ,X_o) =\\
& =  \wh{(Y_1' + Y''_2)}_{[(e, X_o)]_K}\ .
\end{split} \eeq
%
%
%
%
%
%
%
This implies that under the natural isomorphism 
$$\imath: \bR X_o + \gm'  \longrightarrow  T_{y_o} (TS') \ ,$$
 the complex structure  of the  CR structure $(\cD, J)$ of $T G{\cdot} y_o$ preserves the subspace $\gm'$. By counting dimensions  it follows  that  $\gm' $ coincides with  the $J$-invariant subspace  $  \gm$ of $\gl^\perp$ and that 
$ \frac{X_o}{\cB(X_o, X_o)} = \pm Z$.
From this and  \eqref{XdiY}, we also get that $ J \wh Z|_{[(e, X_o)]_K} =  \pi_*\left(0,  \frac{ \pm X_o}{\cB(X_o, X_o)}\right) $. Since the integral curve $\h_t$ of $J \wh Z$  satisfies the conditions  $\h_0 = [(e, 0)]_K$ and $\h_1 = [(e, X_o)]_K$, we get that 
\beq \label{3.6} X_o = Z\ ,\ \   J \wh Z|_{[(e, X_o)]_K} =  \pi_*\left(0, Z\right)\ , \ \ \ \h_t = \left[\left(e,   t Z\right)\right]_K \ .\eeq
\par
\smallskip
\subsubsection{Proof of Theorem \ref{theo32}}
We have now all ingredients for the proof. Consider   the  function  $\t: M_o\setminus S' \longrightarrow \bR$ defined as follows. 
Let  $\h: \bR \to M_o$ be one of the curves defined in \S \ref{optimalcurve} and, for each  $x \in M_o \setminus S'$, let us  denote by $(\gxx, \txx) $ some  pair in  $G \times \bR$  with 
 \beq \label{ecco} x = \gxx {\cdot} \h_{\txx}\ .\eeq
By property (1)  of   $\h$,   such a pair   surely exists, but  is in general not unique. Indeed,  $\gxx$ is  determined  up to composition with some  $ h \in L = G_x$, while  $\txx$ is unique in case   $S'$ is complex and  determined up to a sign in all other cases. 
Then, we set 
\beq\label{deftau}  \t: M_o \setminus S' \longrightarrow (0, + \infty)\ ,\qquad \t(x) \= \left\{\begin{array}{ll} e^{-\txx} & \text{if}\ S'\ \text{is complex}\ ,\\[10pt]
(\txx)^2 & \text{if}\ S'\ \text{is not complex}\ .
\end{array}\right.
\eeq
By construction, $\t$ is $G$-invariant. We claim it  is also   $\cC^\infty$.  By $G$-invariance,    the claim is proven if we show that $\t$ is smooth at each  fixed $y_o \= \h_{t_o}$, $0 < t_o < \infty$.  For this,
consider the decomposition $\gg = \gl + \bR Z + \gm$,  with  $\gl \= \gg_{y_o}$, described in \eqref{decomp}, 
and let
 $\gg^\bC_{y_o}\subset \gg^\bC$ be the  isotropy subalgebra at $y_o$ of the complexified group $G^\bC$.  Further, denote by $\gn$ the  $2n$-dimensional real  subspace $\gn\= \gm + \bR Z + \bR(i Z) \subset \gg^\bC$, which is complementary to $\gg^\bC_{y_o}$, 
and choose  a real basis   $(F_1, \ldots F_{2n})$ 
for $\gn$ with $F_{2n} \= iZ$.  We define $\wt{\exp}_{y_o}: \bR^{2n} \longrightarrow M_o$ by 
\beq
\label{35}  \wt{\exp}_{y_o}(y^1, \ldots, y^{2n}) \=  (e^{y^1F_1} \cdots   e^{y^{2n-1}F_{2n-1}}  e^{y^{2n}F_{2n}}){\cdot} y_o =  (e^{y^1F_1}  \cdots e^{y^{2n-1}F_{2n-1}}) {\cdot} \h_{t_o + y^{2n}}\ .
\eeq
Since $\gn$ is complementary to  $\gg^\bC_{y_o}$,  the Jacobian at $0$ at $\wt{\exp}_{y_o}$  is invertible,  so that, by the Implicit Function Theorem,  the map  $\wt{\exp}_{y_o}$ gives a   diffeomorphism
between a neighborhood   $\cV$ of $0 \in \bR^{2n}$ and a neighborhood $\cU$ of $y_o \in M_o$.  We therefore have that  the inverse
$$\xi = (\wt{\exp}_{y_o})^{-1} : \cU \longrightarrow \cV \subset \bR^{2n}$$
is a    system of real coordinates $\xi= (y^1, \ldots, y^{2n})$ near $y_o$.  
Since $e^{y^1F_1} \ldots e^{y^{2n-1}F_{2n-1}}$ is always an element of the real Lie group $G$, from \eqref{35} we see that  for each  element  $y \in \cU$ we may choose  as $(\gxy, \txy)$ the pair
$$ \gxy =   e^{y^1F_1} \cdots e^{y^{2n-1}F_{2n-1}}\ ,\qquad
  \txy = y^{2n} + t_o\ .
$$
Hence, for each  $y \in \cU$, 
$$\t(y)  = \left\{\begin{array}{ll} e^{-(y^{2n}(y) + t_o)}& \text{if}\ S'\ \text{is complex}\ ,\\[10pt]
(y^{2n}(y) + t_o)^2 & \text{if}\ S'\ \text{is not complex}\ ,
\end{array}\right.$$
from which it follows immediately  that $\t|_{\cU}$ is of class $\cC^\infty$,   as desired. \par
\smallskip
We now want to show that  for any given point  $y \in S'$, there exists a smooth extension of $\t$ on a whole neighborhood of $y$, so that we may consider $\t$ as a smooth real function over  the whole $M_o$. The proof of this property is divided into two cases. \par
\smallskip
\noindent{\it Case 1: the singular $G$-orbit $S'$ is not complex}.   
Under this assumption, either $M_o$ is a Morimoto-Nagano space  or $M  = M_o \cup S$ is an almost homogeneous  manifold with one end and of the mixed type. 
Assume that the first holds, i.e. $ M_o = T S'$ for  a  CROSS of the form  $S' = G/K$. 
By   \S \ref{explicit},  $M_o \simeq G \times_K \gp$ and  $\h_t$ is identifiable with a  curve   \eqref{3.6} for some unitary $Z \in \gp$.  Then,  by $G$-invariance of the norm $\|{\cdot} \| \= \sqrt{- \cB(\cdot, \cdot) }$,    for each 
$x = [(g, Y)]_K \in T S' \setminus\{ \text{zero section}\}$, we may choose as  $(\gxx, \txx)$ the pair
$$\gxx \= g\ ,\qquad \txx \= \pm  \|Y\|^2\ .$$
This implies that the map 
$\t:  M_o \setminus S'  = G \times_K ( \gp \setminus  \{0\}) \to \bR$ has the form
$$\t([(g, Y)]_K) = \|Y\|^2 = - \cB(Y, Y)\ , $$
which  can be directly checked to be $\cC^\infty$ over the entire $T S'  \simeq G \times_K \gp \simeq G^\bC/K^\bC$.\par
\smallskip
 Assume now that  $M = M_o \cup S$ is  an almost homogeneous manifold of the third class. Then  $M_o $ is a $G$-homogeneous bundle over a flag manifold $G/G_Q$, with fibers  given by   Morimoto-Nagano spaces. By $G$-invariance of $\t$ and the fact that curve  $\h$ is entirely contained in a single fiber $\pi^{-1}(z_o)$,  $z_o \in G/G_Q$, the smooth extendibility of $\t$ at the points of the singular $G$-orbit $S'$ is equivalent to the smooth extendibility of the restriction of $\t|_{\pi^{-1}(z_o) \setminus S'}$ to all points  of $\pi^{-1}(z_o)$.  And this   is checked by the same above argument. 
\par
\medskip
\noindent{\it Case 2:  the singular $G$-orbit $S'$ is complex}.  
In this situation,  the  manifold $M$ has two  complex singular $G$-orbits $S$, $S'$ and,  since  it  is in canonical form,  both of them   have real codimension two.
By $G$-invariance of the function $\t$, there is no loss of generality if we  assume that the point $y \in S'$, around which we need to show that $\t$ is smooth,  coincides with the limit point $y = \lim_{t \to + \infty} \h_t$.  Let  $K = G_y$ be the isotropy  at $y$,  so that the singular $G$-orbit $S'$ is identifiable with $S'  = G/K$. Note that, by property (4)  of $\h$ and a dimensional argument, if $\gg = \gl + \gm  + \bR Z$ is the decomposition \eqref{decomp} corresponding to the regular point $x_o = \h_0 \in M_o \setminus S'$, then
  $\gk  = \gl + \bR Z$.
  \par
  Consider  a $K$-invariant K\"ahler metric  $g$ on $M$ and  let   $V = (T_y S')^\perp$ be the $2$-dimensional $g$-orthogonal  complement to $T_y S'$  in $T_y M_o$. Since $S'$ is complex, $V$ is $J$-invariant.  Denote by  $\exp_y: V \to M$ 
  the restriction to $V $ of the exponential map of $(M, g)$ at $y$. \par
  \smallskip
    By standard facts  on  proper actions (see  e.g. \cite{Br, GO, PoS0}),   the isotropy representation of  $K$ on $T_y M$ preserves $V$  and acts  linearly and   isometrically on  $V$ with codimension one regular orbits. It also preserves the complex structure $J_o = J_y|_V$. Hence,  we may  identify $(V, \Jst, g_y|_V)$ with $(\bC, \Jst, \langle \cdot, \cdot \rangle)$, where $ \langle \cdot, \cdot \rangle$ is  the standard Euclidean product, so that  the representation of $\rho: K \to \GL( V)$ is such that $\rho(K) = S^1$. 
\par
  \smallskip
  Consider  the linear bundle $\pi: G \times_K \bC \to S' = G/K$. It is known that   there exists a $K$-invariant neighborhood $\cU \subset V$ of the origin such that the map 
\beq \varphi: G \times_K \cU \longrightarrow M_o\ ,\qquad \varphi([g, v]_K) \=  \exp_{g {\cdot} y}(v) = g {\cdot} \exp_y(v)\eeq
is a $G$-equivariant diffeomorphism between $G \times_K \cU$ and a neighborhood $\cV = \varphi(G \times_K \cU)$ of $S'$ in $M_o$. By construction, the singular orbit $S' $ coincides with the image by $\f$ of the zero section $G\times_K \{0\}$ and the action of   $e^{\bR Z}$ on $\varphi(\cU)$ corresponds to the action on $G \times_K \cU$ defined by  
$$e^{tZ}{\cdot} [(g,\z)]_K\=  [(g,e^{it} \z)]_K\ .$$
This yields that  the infinitesimal transformations  $\wh  Z$,  $J \wh Z = \wh{(i Z)}$ determine  on each fiber $\{g\} \times_K \cU \simeq\cU \subset  \bC$ of $G \times_K \cU$  the  (real) vector fields  of $\bC \simeq \bR^2$
$$\wh Z = i \z \frac{\p}{\p \z} -i  \overline \z  \frac{\p}{\p \bar  \z} \ ,\qquad J \wh  Z = - \z \frac{\p}{\p \z} +  \overline \z  \frac{\p}{\p \bar  \z}  \ .$$
The  flow of $J \wh  Z $ in $\{g\} \times_K \cU \simeq\cU \subset  \bC$  is  then given by  
 $$\Phi^{J \wh Z}_t(\z) = e^{ - t} |\z|^2\ .$$
  From this, we see that $\h$ corresponds to the curve on $G \times_K \cU$
  $$\check\h_t\= \varphi^{-1}(\h_t) = [(e, e^{-t} \z_o)]_K \qquad \text{where}\ \z_o\ \text{is such that}\ \ [(e, \z_o)]_K = \f^{-1}(x_o)\ .$$
So,  for each $x =  \varphi([(g, \z)]_K)$, $\z \neq 0$, in $\f(\cU)$, we may choose    $(\gxx, \txx)$ as the pair
 $$\gxx = g\ ,\qquad \txx(\z) = -\frac{1}{2} \log\left(\frac{ |\z|^2}{|\z_o|^2}\right)\ . $$
Hence,  for all  points of $\cV\setminus S' \simeq \varphi(G \times_K (\cU \setminus \{0\}))$,  the function $\t$ is such that
$$ \t(\varphi([(g,\z)]_K)  = C |\z|^2\qquad \text{with}\ \ C \= \frac{ e^{\frac{1}{2}}}{|\z_o|^2}$$
and it  is   clearly smoothly extendible  at all points of $ \cV \simeq  \varphi(G \times_K \cU)$, as claimed. \par
\smallskip
We now want  to show that, in all cases,  the map  $\t: M_o \to [0, + \infty)$ satisfies (1) and (2) of 
Definition \ref{defMonge}. For  this, consider the distributions $\cZ$, $\cH \subset  T(M_o\setminus S')$,  defined by 
$$\cZ_x \= \big\langle \wh Z_x, J \wh Z_x\big\rangle \ ,\qquad \cH_x \= \cD_x \ ,\qquad x \in M_o \setminus S'\ ,$$
where  $Z$ is the  unitary element  of $\gg$,   appearing in  \eqref{decomp}, which is  $\cB$-orthogonal to the isotropy   $\gl = \gl_x$ and to the space $ \gm \simeq \cD_x$, corresponding to  the CR structure of  the orbit $G{\cdot} x$. Note that such  $Z \in \gg$ {\it does depend on   $x$} and, for the sake of  clarity,  it will be later denoted  by $\Zxx$. \par
A direct check shows that  $\cZ$ and $\cH$ are $d d^c \t$-orthogonal and that, for each $x$,   the restricted $2$-form  $ d d^c \t_x|_{\cH_x\times \cH_x}$ 
 is (up to a multiple) the Levi form of the $G$-orbits $G{\cdot}x$, hence    strictly positive. So,  for (1) of Definition \ref{defMonge},  we just need to show that  
 $ d d^c \t|_{\cZ_x\times \cZ_x}$ is positive at each $x \in M_o \setminus S'$. For this, we first observe that   the distribution $\cZ$ is integrable and its integral  leaves are the orbits in $M_o \setminus S'$ of the complex Lie groups $\exp(\bC \Zxx)$.  Hence
 \beq \label{eccoqua0}  d d^c \t|_{\cZ_x\times \cZ_x} = d d^c (\t|_{\exp(\bC \Zxx){\cdot}x})\big|_x\ .\eeq
By previous discussion,   
\beq
 \label{eccoqua} \t|_{\exp( \z \Zxx){\cdot}x} = C |\z|^2\qquad \text{or} \qquad  \t|_{\exp( \z \Zxx){\cdot}x} = C (\Im \z)^2\ ,
\eeq
the first occurring when $M = M_o \cup S$ has two complex singular $G$-orbits, the second in all other cases.
From \eqref{eccoqua0} and  \eqref{eccoqua},  we get  $ d d^c \t_x|_{\cZ_x\times \cZ_x}> 0$ in all cases,  as desired.\par
In order to check  (2) of Definition \ref{defMonge}, we first observe that  for any smooth  function $f: (0, + \infty) \to \bR$ with nowhere vanishing differential $df$, we have 
$d d^c (f \circ \t)|_{\cH_x\times \cH_x} > 0$ at each point $x$.
Therefore,  by \eqref{eccoqua0},   property (2) holds if and only if there is such  an $f$   with  harmonic restrictions $f \circ \t|_{\exp(\bC \Zxx)}$.  By \eqref{eccoqua}, we see that   $f(t) = \log(t)$ and $f(t) = \sqrt{t}$ satisfy the request in the two cases. This concludes the proof  of (i). \par
\smallskip
 It remains to prove (ii). Let $\t': M_o \to [0, + \infty)$ be a $G$-invariant Monge-Amp\`ere exhaustion of $M_o$ satisfying ($\a$), ($\b$) and ($\g$). Since $d d^c \t'$ and $d d^c (\log \circ \t')$ are both  positive on the CR distributions of  the  regular $G$-orbits (which are  level sets of $\t'$ and $f \circ \t'$ and are strongly pseudoconvex),  due to ($\alpha$)  the restrictions  $\log \circ \t'|_{\exp(\bC \Zxx){\cdot}x}$ or $\sqrt{\t' }|_{\exp(\bC \Zxx){\cdot}x}$ are  necessarily harmonic. They are also  constant along   the sets  $\exp( \bR \Zxx){\cdot}x $, which are the intersections of $(\exp(\bC \Zxx) {\cdot} x)$ with the $G$-orbits.\par
Let us   first prove that all this implies (ii)  when $M = M_o \cup S$ has two ends.  By the proof of (i), we know that  each orbit $\exp(\bC \Zxx){\cdot}x$ is identifiable with  $\bC \setminus \{0\}$, so  that  the group $\{\exp(\z \Zxx), \z \in \bC\}$ corresponds to  the group  
 $\{ D_\z(z) \= e^{i \z}{\cdot} z, \ \z \in \bC\}$.  Under such identification, the restriction $\log \circ \t'|_{\exp(\bC \Zxx){\cdot}x}$  is an harmonic function of $\bC\setminus \{0\}$  depending only on the distance from the origin. It has therefore the form
 \begin{align*}
 & \log\circ \t'(\exp(\z \Zxx){\cdot}x) = a_x  + 2 b_x \log(|\z|)\qquad \text{for some constant}\ a_x, b_x \in \bC\ ,\ \text{so that}\\
&  \t'(\exp(\z \Zxx){\cdot}x) = c_x (|\z|^{2})^{b_x}\qquad \text{with}\ c_x \= e^{a_x}\ .%
\end{align*}
Condition  ($\g$) and $G$-invariance imply  that   $b_x \equiv  1$ and that the constant  $c_x$ does not depend on $x$.  From   \eqref{eccoqua},  (ii)  follows in this  case.\par
We now prove   (ii) holds when $M$ has only one end. In this case,  each orbit $\exp(\bC \Zxx){\cdot}x$ is  identifiable with  a quotient $\bC/\G$ with   $\G$ group of real translations $\G =\{T_k(z)\= z + 2 \pi k \z\ ,\ k  \in \bZ\}$. Under this identification,  the group $\{\exp(\z  \Zxx)$ acts on such orbit as the group of complex translations 
$\{ T(z)\= z + \z\ ,\ \z \in \bC\}$ and the restriction   $\sqrt{ \t'}|_{\exp(\bC \Zxx){\cdot}x}$  is identified with  an harmonic function of $\bC/\G$, which is constant on the lines $\{\Im(z) = c\}$, $c \in \bR$.  Hence
$$  \sqrt{ \t'(\exp(\z \Zxx){\cdot}x)} = a_x  + 2 b_x \Im(\z)\qquad \text{for some constant}\ a_x, b_x \in \bC\ .$$
Condition  ($\g$) and $G$-invariance imply  that   $a_x \equiv  0$ and that the constant  $b_x$ does not depend on $x$. From    \eqref{eccoqua},    claim (ii)  follows also in this case. 
\par
\medskip
\section{Deformability versus Rigidity}
\setcounter{equation}{0}
As we already mentioned, Theorem \ref{theo32} is a source of   several  new examples of Monge-Amp\`ere spaces. On the other hand, 
 the   previously known examples  include   two important families of Monge-Amp\`ere spaces,    the manifolds of circular type
and the Grauert tubes,  having the following contrasting  properties:  the first can be 
all considered  as deformations of   $\bC^n$ or  $\bB^n$ (\cite{PS, PS2}),   while the second are  characterized by   strong rigidity results  (\cite{LS}). Motivated by this,   in this final  section  we  investigate  whether also the new examples enjoy   manifest   deformability (or rigidity) properties.
Let us begin  by fixing      the meaning of  ``deformability''  for a Monge-Amp\`ere   space.  \par
\smallskip
\subsection{Riemann mappings and deformations of modeling spaces}
\label{Riemannmappings}
\hfill\par
Let $\cX$ be a Monge-Amp\`ere space,   which is a Remmert reduction  $\pi: \wt \cX \to \cX$ of a complex manifold $\wt \cX$  with Monge-Amp\`ere exhaustion $\t: \wt \cX \to [0, T)$. Let also $u = f \circ \t: \wt S \to \bR$ be  function satisfying the conditions $2 i \partial \overline \partial  u  \geq 0$ and $(\p \overline \p  u)^n = 0$,  and  $ \cS = \pi(\t^{-1}(0)) $  the soul of $\cX$. We recall that $\cX \setminus \cS$ is a complex manifold, naturally identifiable with  $\wt \cX \setminus \t^{-1}(0)$, so that  both exhaustions  $\t$ and $u$ are well defined and smooth on $\cX \setminus \cS$.\par\smallskip 
We now observe that all {\it local} properties of  exhaustions $\t$ and  $u$  on the manifold $\cX \setminus \cS$, which  have been proven  in the literature for some specific cases (as, for instance, when  $\cX$ is a  manifold of circular type or a Morimoto-Nagano spaces  -- see e.g. \cite{St1, Pt1, PS1, PW}) are 
valid for  {\it any}  Monge-Amp\`ere space $\cX$. In particular,  one can directly check that   there is always a well defined  vector field $Z$  on $\cX \setminus \cS$ that satisfies the condition 
 \beq \label{21} d d^c \tau(J Z, J X) =  X(\tau) \ \qquad \text{for any vector field}\ \ X \in T\left( \cX \setminus \cS\right)\ .\eeq
For the spaces  of  Theorem \ref{theo32},  this vector field  verifies  $Z_x = \pm \wh{Z(x)}_x$ at each $x \in M_o \setminus S'$. 
By  $J$-invariance of the $2$-form $d d^c \t$ and integrability of the complex structure $J$,  the vector field $ Z$ is tangent to each level set $\t^{-1}(c)$, $c \in (0, T)$, and  $Z$ and  $J Z$ generate a $J$-invariant, integrable $2$-dimensional distribution $\cZ \subset T (\cX \setminus \cS)$, called {\it Monge-Amp\`ere distribution}. A direct computation shows that it  coincides with the distribution defined by $ \cZ_x = \ker \left.d d^c u\right|_x$ at each $x \in \cX \setminus \cS$.\par
The   foliation $\cF$, given  by  the integral leaves of $\cZ$,  is called  the {\it Monge-Amp\`ere foliation  of $(\cX, \t)$}.  The closures in $\cX$ of these leaves sometimes form  a regular foliation,  sometimes  not.  However, in all known examples,  the closures of their   lifts on the manifold   $\wt \cX$   always form  a regular foliation of $\wt \cX$. For brevity,   when this property occurs we say that   the Monge-Amp\`ere foliation is {\it $\wt \cX$-regular}. 
For the spaces  in  Theorem \ref{theo32}, the Monge-Amp\`ere foliation $\cF$  consists of  orbits in $M_o\setminus S'$ of the  $1$-dimensional complex groups $\exp(\bC \Zxx)$ and, in all such cases,  $\cF$   is $M_o$-regular.\par
\smallskip 
It is clear that  if  $\f: \cX \to \cX'$ is a biholomorphism  between two  Monge-Amp\`ere spaces $(\cX,   \t)$,  $(\cX',  \t')$ with   $\t = \t' \circ \f $,  then $\f$  maps biholomorphically  each  leaf of the Monge-Amp\`ere foliation foliation of  $\cX$ into  a corresponding leaf of  the Monge-Amp\`ere foliation of $\cX'$ . 
In certain cases this property admits an inverse, in the sense  that   if  $\f: \cX \to \cX'$ is an homeomorphism with  $\t = \t' \circ  \f$ and  mapping biholomorphically  each Monge-Amp\`ere leaf of $\cX$ into a corresponding leaf of $\cX'$, then $\f$ is  a  biholomorphism provided that certain additional hypothesis are satisfied. Stoll's characterization of $\bC^n$ and 
 Lempert and Sz\"oke rigidity theorems for Grauert tubes  can be considered as examples  of  such kind of property.  All this motivates   the next notion. \par
\smallskip
 Consider   a fixed  Monge-Amp\`ere space  $(\cX_o,  \t_o)$, which we call {\it model}  from now on.  Let also  $\pi: \wt \cX_o \to \cX_o$ be the  Remmert reduction that determines the Monge-Amp\`ere space $\cX_o$ and    denote  by  $\cS_o $  and $\cF_o$ the soul  and   the Monge-Amp\`ere foliation of $\cX_o$, respectively.    We assume that $\cF_o$ is $\wt \cX_o$-regular, as it occurs in all considered examples.\par
\begin{definition}  Let  $(\cX,  \t)$ be a Monge-Amp\`ere space,   Remmert reduction of a complex manifold $\wt \cX$,   with soul $\cS$ and  $\wt \cX$-regular Monge-Amp\`ere foliation $\cF$.
 We say that $(\cX,  \t)$ is {\it modeled on  $(\cX_o,  \t_o)$} if there is a homeomorphism $\f: \cX_o \to \cX$ such that  $\t_o = \t \circ \f$ and: 
 \begin{itemize}
 \item[i)]   $\f|_{\cX_o \setminus \cS_o}: \cX_o \setminus \cS_o \to \cX \setminus \cS$ is a diffeomorphism   mapping biholomorphically  each  leaf of   $\cF_o$ into  a  leaf of $\cF$; 
 \item[ii)] the  restriction  $\f|_{\cX_o \setminus \cS_o}$  lifts to a diffeomorphism 
 $\wt \f: \wt \cX_o \setminus \pi^{-1}(\cS_o) \to \wt \cX \setminus \pi^{-1}(\cS)$
 which  smoothly extends to a diffeomorphism  between the $\wt \cX_o$ and $\wt \cX$.
 \end{itemize}
   Any such homeomorphism $\f: \cX_o \to \cX$ is called {\it Riemann mapping}  of $\cX$.  The spaces that are modeled on  $(\cX_o,  \t_o)$, but are not biholomorphic to $\cX_o$, are called  {\it non-trivial deformations} of the model. 
\end{definition}
\par
\smallskip

\subsection{Soul rigidity, soul semi-rigidity  and free deformability}\hfill\par
Looking at the known examples, there are  models with a lot of non-trivial deformations and others with no deformation with sufficiently regular Riemann mappings. For instance,  the results in \cite{BD, Le, PS}    show that any  smoothly bounded strictly convex  domain in $\bC^n$,   not biholomorphic to $\bB^n$,  is    a non-trivial deformation of the standard unit ball $(\bB^n, \Jst, \|{\cdot}\|^2)$.  On the other hand, Stoll's characterization of $\bC^n$   (\cite{St1, Bu, PS2})  shows 
 that there exists no non-trivial deformations  of  $(\bC^n, \Jst, \|{\cdot}\|^2)$ in  the class of  Monge-Amp\`ere manifolds  with  Riemann mappings of class   $\cC^2$.  Known rigidity results  for Grauert tubes give uniqueness for the Riemann mappings from  Morimoto-Nagano spaces (\cite{GS, LS, PW}).  For  clarifying similarities and differences between  all such results, we now introduce the following 
 notions.  \par
 \smallskip
 Let $\f: \cX_o \to \cX$ be a Riemann mapping from a model $\cX_o$  and denote by $\wt \f: \wt \cX_o \to \wt \cX$ the corresponding 
 lifted diffeomorphism between the manifolds,   which project onto  the two spaces  by  Remmert reductions $\pi_o: \wt \cX_o \to \cX_o$ and $\pi: \wt \cX \to \cX$. Denoting by  $J_o$, $J$  the complex structures of $\wt \cX_o$, $\wt \cX$,  respectively,  the diffeomorphism $\wt \f$ is a biholomorphism if and only if for any tangent vector $v \in T_x \wt \cX_o$, $x \in \wt \cX_o$, 
 \beq \label{biho} \f_*(J_o(v)) = J \f_*(v)\ .\eeq
It is now convenient  to  consider the following  weaker conditions. We say that  $\f$ is:
 \begin{itemize}
 \item[1)] a {\it biholomorphism at the blow ups of the souls} if it satisfies \eqref{biho} for  any  vector $v  \in T_x \pi^{-1}_o(\cS_o)$  in a tangent space of the  preimage $\pi^{-1}(\cS_o)$ of the soul $\cS_o$;
 \item[2)] a {\it biholomorphism between the souls} if it satisfies \eqref{biho} for  any  vector $v  \in T_x \pi^{-1}_o(\cS_o)$  in a tangent space of  $\pi^{-1}(\cS_o)$ which projects onto a non trivial tangent vector  of   $\cS_o$.
 \end{itemize}
The second condition is manifestly weaker than the first since it does not requires that \eqref{biho}  holds for  vectors that are in  $\ker \pi_*|_{y}$, $y \in \pi^{-1}(\cS_o)$. Moreover, note that: 
\begin{itemize}
\item[a)] When $\cX_o $ is the  Remmert reductions of  some  $M_o = M \setminus S$, with $M$ almost homogeneous  with two ends, condition (2) is trivial, since in this case   the soul   is an isolated point and there are no non-trivial tangent vector for such soul.
\item[b)]ÊWhen $\cX_o$,  $\cX$ are Grauert tubes over two Riemannian manifolds $(M_o, g_o)$,  $(M, g)$, respectively, condition (2) coincides with condition  (1) and  it is equivalent to require that the Riemann mapping $\f$ induces an isometry between $(M_o, g_o)$ and $(M, g)$. 
\item[c)]ÊWhen $\cX_o$ and $\cX$ are  manifolds of circular type, as e.g. two strictly convex domains in $\bC^n$, (2) is a trivial condition because of (a), while  (1) is equivalent to require that the Riemann mapping induces an isometry between the Kobayashi indicatrices at the  centers. 
\end{itemize}
We now say that a model $\cX_o$ is  \par
 \begin{itemize}
 \item {\it soul rigid} if any  Riemann mapping $\f: \cX_o \to \cX$ which is a biholomorphism between  the souls is a biholomorphism between the two Monge-Amp\`ere spaces;
\item {\it soul semi-rigid}  if it is not soul rigid, but nonetheless for any Riemann mapping  $\f: \cX_o \to \cX$ which is a biholomorphism at the blow ups of the souls is a biholomorphism between the two Monge-Amp\`ere spaces;
 \item {\it fully deformable}  if  it is  not of the previous two types. 
 \end{itemize}
\par
The quoted rigidity results for Morimoto-Nagano spaces and Grauert tubes can be stated  saying that  those  manifolds are soul rigid.  On the other hand, 
 the results of   \cite{Bu, Pt1} show that  
 $(\bC^n,  \|{\cdot}\|^2)$ is  soul semi-rigid, while    the examples in  \cite{BD, PS, PS2} show  that the standard unit ball $(\bB^n, \Jst, \|{\cdot}\|^2)$ is a fully deformable model with a lot of non-trivial deformations. \par
\smallskip
\subsection{New examples of semi-rigid and fully deformable  models}\hfill \par
The following theorem   gives a  common framework for the  so far known  rigidity results on Monge-Amp\`ere spaces. It also  indicates a new interesting class of   semi-rigid examples and,  by   the examples  in \S \ref{counterexamples} below,  suggests  the existence of a large new family of   fully deformable  Monge-Amp\`ere spaces.\par
\begin{theo} \label{main}Let   $M$ be  one of the almost homogeneous manifolds  in canonical form,   described  in \S\S\,\ref{3.2.1}, \ref{3.2.2},  \ref{3.2.3},  and   $\t_o: M_o = M \setminus S \to [0, + \infty)$ the Monge-Amp\`ere exhaustion  of Theorem \ref{theo32}. 
Let also $\cM_o$  be the corresponding Remmert reduction of $M_o$,  equipped  with  the  exhaustion induced by $M_o$,  which,  for simplicity of notation, we also  denote by $\t_o$.
For any $c \in (0, + \infty]$, let $\cM_o{(c)}\= \{ x :  \t_o(x) < c\}$ be  the Monge-Amp\`ere subspace   with exhaustion $ \t_o|_{\cM_o(c)}$. Then: 
\begin{itemize}
\item[i)] if $M$ has two ends,  the space $ (\cM_o{(\infty)}, \t_o) = (\cM_o, \t_o) $  is soul semi-rigid;  on the other hand,  there is an  $M$ with two ends for which  all   Monge-Amp\`ere spaces  $(\cM_o(c), \t_o)$  with $c < \infty$  are fully deformable; 
\item[ii)] if $M$ is a compactification of a Morimoto-Nagano space, each  Monge-Amp\`ere space  $(\cM_o(c), \t_o)$,   $0 < c \leq \infty$,  is soul rigid; 
\item[iii)] if $M$ has one end and   is  of mixed type, then  each Monge-Amp\`ere space  $(\cM_o(c), \t_o)$, $0 < c \leq +\infty$,  is soul  semi-rigid. 
\end{itemize}
\end{theo}
The proof is crucially based on some properties   of  ``deformation tensors''  of certain complex structures and on some known  counterexamples to soul rigidity or semi-rigidity.   Before  proving this theorem, we need to review   such results in some detail. \par
\smallskip
\subsubsection{Deformation tensors of deformed  Monge-Amp\`ere spaces}\label{deformationtensors}
 As usual,  let $\cX$ be a Monge-Amp\`ere space, which is  Remmert reduction   of a complex manifold $\wt \cX$  with   Monge-Amp\`ere exhaustion $\t: \wt \cX \to [0, T)$. Denoting by $\cS$ the soul of $\cX$, we observe that   on $\cX \setminus \cS $ (which always   identifiable with the complementary set  $\wt \cX \setminus \pi^{-1}(\cS)$)    one can consider  the  $J$-invariant  distribution $\cH$, 
called  {\it normal  distribution\/},  
 defined by 
\beq \cH_x = \{\ X \in T_xM\ :\ dd^c\tau(Z, X)_x = dd^c\tau(J Z, X)_x = 0 \ \}\ .\label{2.12}\eeq
By  non-degeneracy of $d d^c \t$ on $\cX \setminus \cS$,  we have  that $T_x M = \cZ_x \oplus \cH_x$ at each $x \in \cX \setminus \cS$. Moreover,  for each level set  $ \t^{-1}(c)$, $c \in (0, T)$, the restriction $\cH|_{ \t^{-1}(c)}$  coincides with  the $J$-invariant distribution of  the induced CR structure $(\cD = \cH|_{ \t^{-1}(c)}, J)$ of such level set.  
\par
\smallskip 
Consider now  a Monge-Amp\`ere space   $(\cX, \t)$ modeled on a manifolds $ (\cX_o,  \t_o)$,  and denote by   $J_o$,  $J$  the complex structures of the complex manifolds $\wt \cX_o$,  $\wt \cX$,  of which $\cX_o$ and  $\cX$ are Remmert reductions, respectively.  Fix also a Riemann mapping
$\varphi: \cX_o\to \cX$, with associated  lifted diffeomorphism  $\wt \varphi:  \wt \cX_o  \to \wt \cX$. As direct consequence of definitions, $\wt \f$ sends the Monge-Amp\`ere and normal distributions of $\cX_o \setminus \cS_o$ into the corresponding distributions of $\cX \setminus \cS$. \par
\smallskip
We denote by $\bJ  \=  \wt \f^{-1}_*(J)$ the pull-back  on $\cX_o$ of the complex structure $J$. Clearly  $\f$ is a biholomorphism if and only if $J_o = \bJ$.
 \par
\smallskip
We  recall that,  being a Riemann mapping, the map $\wt \varphi$ is a biholomorphism along each leaf of the Monge-Amp\`ere foliation. This mean that   $J_o|_{\cZ} = \bJ|_{\cZ}$ and  that 
differences between  $J_o$ and $\bJ$  might occur only when they are restricted on the normal distribution $\cH$. On the other hand, at each point $x$, 
both complex structures  $J_{o}|_{\cH_x}, \bJ|_{\cH_x}: \cH_x \to \cH_x$  are  determined by their $(-i)$-eigenspaces in the complexification $\cH^\bC_x$.  Let us  denote these eigenspaces by $\cH^{01}_x$ and  $\cH'{}^{01}_x$, respectively, and indicate by $\cH^{01}, \cH'{}^{01} \subset \cH^\bC$ the  involutive complex distributions, determined by such  $(-i)$-eigenspaces.  Their conjugate distributions (which are generated by the $(+i)$-eigenspaces)  are  denoted by $\cH^{10} = \overline{\cH^{01}}$   and  $\cH'{}^{10} = \overline{\cH'{}^{01}}$. \par
We now  observe  that,  for a fixed  $x$,  if the standard projection $p: \cH^\bC_x = \cH^{10}_x+ \cH^{01}_x  \to \cH^{01}_x$ satisfies the condition  $p(\cH'{}^{01}_x) =  \cH^{01}_x$, then the   space $\cH'{}^{01}_x$ has the form
\beq \label{def1}\cH'{}^{01}_x = \left\{v = w+ \phi_x(w)\ ,\ w \in \cH^{01}_x\right\}\eeq
for some appropriate  tensor
$\phi_x \in (\cH^{01}_x)^* \otimes \cH^{10}_x $. The set $\cU \subset \cX(c)$ of  points $x$,  for which the condition  $p(\cH'{}^{01}_x) =  \cH^{01}_x$ holds,  is  open  and we call it  the {\it regularity set of $\bJ$}. The   tensors $\phi_x$   combine to a smooth tensor field  $\phi$ on  $\cU$,  called  the  {\it deformation tensor of $\bJ$} (\cite{BD, PS, PS2}).   \par
The  biholomorphicity condition can be now expressed in terms of the deformation tensor $\phi$   saying that  {\it   $\f$ is a biholomorphism if and only if  the regularity set $\cU$ of $\bJ$   coincides with   $\wt \cX_o$ and  
$\phi \equiv 0$}.  \par
 \smallskip
 We conclude this section recalling a crucial property of  the distributions $\cZ$ and $\cH$.
  We  remark that the normal distribution $\cH$ is  invariant under the 
 flows of the vector fields $Z$ and $JZ (= \bJ Z)$, defined in \eqref{21}.  Due to  this and the fact that $\cZ^{01} + \cH'{}^{01}$ is involutive (due to the integrability of $\bJ$),  for each vector field in $Y \in \cH^{01}|_{\cU}$
\beq   [Z^{01}, Y + \phi(Y)]  = [Z^{01}, Y] + \phi([Z^{01}, Y] )\ , \qquad \text{where}\ Z^{01} \= Z + i J Z\ .\label{5.9}
\eeq
So, if we denote by  $L$ an  integral complex leaf  of  $\cZ$ and by $(e_\a, e_{\bar \b} \= \overline{e_\b})$ is    a   frame field  for $\cH^{\bC} = \cH^{10}+ \cH^{01}$ on a neighborhood of $L$,   invariant under the complex flows of  $Z$, and if  $\z$ is a complex coordinate on $L$  with $\frac{\p}{\p \z} = (Z - i JZ) = \overline{Z^{01}}$, then \eqref{5.9} yields that  the components  $\phi_{\bar \a}^\b|_L$ of    $\phi|_L$ in the frame $(e_\a, e_{\bar \b})$ are   holomorphic functions of $\z$, i.e. 
$$ \frac{\p \phi_{\bar \a}^\b}{\p \overline \z}\bigg|_x = 0\qquad \text{at each}\ x \in L\ .$$
\par
\smallskip
\subsubsection{Counterexamples to soul rigidity and soul semi-rigidity} \label{counterexamples}Assume that   $M$ is an almost homogeneous manifold as in (i) or (iii) of Theorem \ref{main}.   Then $M$ is a  fiber bundle $p: M \to G/K$ over a flag manifold $(G/K, J_{G/K})$, with  fiber $F$ equal   either  to $\bC P^1$ or to a compactification of  a Morimoto-Nagano space. 
In all these cases,  the projection $p$ is holomorphic with respect to the $G$-invariant complex structures  $J$ of $M$ and 
$J_{G/K}$ of $G/K$. Now, 
given  an open set $\cU \subset G/K$, for which  there exists  a holomorphic trivialization $p^{-1}(\cU) \simeq \cU \times F$,   we may consider a non-trivial non-biholomophic diffeomorphism 
$h : G/K {\to} G/K$  mapping $\cU$ into itself and satisfying  the  triviality condition  $h|_{G/K \setminus \cU}= \Id|_{G/K \setminus \cU}$ on the complement of $\cU$. 
Since $p^{-1}(\cU)$ is holomorphically trivializable, we may construct  a fiber preserving diffeomorphism $\wt \f: M \to M$ such that: 
\begin{itemize}
\item[A)] it projects onto 
$h$ and is  such that $\wt \f|_{M \setminus p^{-1}(\cU)} = \Id_{M \setminus p^{-1}(\cU)}$; 
\item[B)] it is holomorphic on each fibre of $p: M \to G/K$ and maps the level sets $\t^{-1}(c) \cap F_x$ of each fiber $F_x = p^{-1}(x)$, $x \in \cU$, into the corresponding level sets $\t^{-1}(c) \cap F_{h(x)}$ of the fiber $F_{h(x)} = \pi^{-1}(h(x))$. 
\end{itemize}
Since all  fibers $F_x$, $x \in \cU$,  are  identifiable  one to the other  by means of  a holomorphic  trivialization $\pi^{-1}(\cU) \simeq \cU \times F$, a  map that satisfies (A) and (B) can be easily determined.\par 
By construction, such diffeomorphism $\wt \f: M \to M$ leaves invariant all submanifolds $M_o(c) \= \pi^{-1}(\cM_o{(c)})$, $c \in (0, + \infty]$,  and  its restrictions to them  are all Riemann mappings. But no such restriction is  a biholomorphism. On the other hand, by construction, it is  a biholomorphism between the souls,  since it satisfies \eqref{biho} for any vector $v \in T_{\pi^{-1}(x)}(\pi^{-1}(\cS_o))$ which  projects onto some non trivial tangent vectors of the soul $\cS_o$, this being  diffeomorphic to each of the zero level sets  $\t^{-1}(0) \cap F_x$ of all fibers of the fibration $p: M \to G/K$ 
\par
\smallskip
Assume now that $M$ is the blow up $p: \wt{\bC P^n} \to \bC P^n$ of $\bC P^n$ at some point $x_o$.  As we already mentioned,  this is an almost homogeneous manifold with two ends, acted on by $G = \SU_n$ and such that the manifold $M_o = M \setminus S$ is naturally identifiable with the blow up  $p: \wt{\bC^n} \to \bC^n$ of $\bC^n$ at the origin. In this case, the Monge-Amp\`ere spaces $(\cM_o{(c)}, \t)$, $0 < c < + \infty$,  are identifiable with  the  balls $\bB^n_c \subset \bC^n$ of radius $c$ and center at  $0$, equipped with the standard exhaustion $\|{\cdot}\|^2$. By \cite{PS2}, Cor. 6.2 and Thm. 5.2, we know that there are Monge-Amp\`ere spaces,   biholomorphic to  strictly convex domains, which are  non-trivial deformations of $(\bB^n_c, \|{\cdot}\|^2)$ with  Riemann mappings that  induce the identity map on the pre-images of the souls. This means that all such Monge-Amp\`ere spaces are fully deformable models.  We  expect that  this property holds for {\it all} models $(\cM_o{(c)}, \t)$, $0 < c < + \infty$,  that are determined by the  almost homogeneous manifolds with two ends in canonical form.  \par
\smallskip
\subsubsection{Proof of Theorem \ref{main}.} Let  $\f: \cM_o{(c)} \to \cX$  be a  Riemann mapping between one of the  models $(\cM_o{(c)}, \t_o)$  considered in (i) - (iii) and a Monge-Amp\`ere space $(\cX,  \t)$. As usual, we denote by   $\wt \f: M_o{(c)}\to \wt \cX$  the associated lifted diffeomorphism and by $\bJ = \wt \f^{-1}_*(J)$ the pull-back on $M_o{(c)} $ of the complex structure $J$ of $\wt \cX$.
 We first want  to prove that, for all models in (ii) and (iii),   if $\f$  satisfies \eqref{biho} at all tangent spaces of $\pi^{-1}(\cS_o)$, then  
 the regular set $\cU$ of $\bJ$ coincides with  the whole $M_o{(c)}$ and the corresponding deformation tensor $\phi $ vanishes identically, meaning that  $\f$ is a biholomorphism.\par
 \smallskip
We recall  that  the distributions $\cZ$ and $\cH$ are  both $J_o$ and $\bJ$-invariant and that the complex structures $J_o$ and $\bJ$ agree on the vector fields in $\cZ$. Note also that for each   $ \cM_o(c)$,  the  $J_o$-invariant  distributions $\cZ$ and $\cH$, taken as distributions on $M_o(c) \setminus \pi^{-1}(\cS) = M_o(c) \setminus S'$,  extend smoothly at all points of  $M_o(c)$.  
These property imply that $\f$  satisfies \eqref{biho} at all tangent spaces of $\pi^{-1}(\cS_o)$ if and only if  $J_o|_{\cH_x} = \bJ|_{\cH_x}$ at all points $x \in \pi^{-1}(\cS_o) = S'$. 
If this is the case, then the regular set $ \cU$ of the complex structure $\bJ$ clearly includes the   submanifold $\pi^{-1}(\cS_o) = S'$,   the  deformation tensor $\phi$ of $\bJ$ is well defined on  a tubular neighborhood $\cW$ of $\pi^{-1}(\cS_o)$ and the restriction  $\phi|_{\pi^{-1}(\cS_o)}$ is identically equal to $0$. \par
\smallskip
Assume now that $M = M_o \cup S$ is either a compactification of a Morimoto-Nagano manifold or an almost homogeneous manifold with one end and of mixed type. In these two cases, the leaves of the Monge-Amp\`ere foliation of $M_o{(c)} \subset M_o$ are (contained in) orbits of the $1$-dimensional Lie groups $\exp(\bC \Zxx)$ described in \S \ref{Riemannmappings}, and  intersect the pre-image $\pi^{-1}(\cS_o) = S'$ of the soul  along sets, which have  Hausdorff dimension $1$. This fact, together with 
 the holomorphicity of the components of $\phi$ in the complex  coordinate 
$\z$ of each leaf $L = \{\exp(\z \Zxx){\cdot} x, \z \in \bC\}$,  implies that if  the Riemann mapping $\f$  is  a  biholomorphism at the blow ups of the souls,  then the restriction  $\phi|_L$ is identically equal to $0$ along each such leaf $L$. This means that the regular set $\cU$ of $\bJ$ contains all leaves of the Monge-Amp\`ere foliation, hence the whole  $\cU =M_o{(c)}$,   and that $\phi$ vanishes identically  on  $M_o{(c)}$ as claimed. \par
\smallskip
We  claim that the same conclusion holds also if   $M$ is as in (i) and the considered model is   $(M_o{(c)}, \t)$ with $c = \infty$. In this case the  closure   in $M_o{(\infty)} =  M_o$ of a leaf $L$   intersects  the pre-image $\pi^{-1}(\cS_o)$ in a single  point. So, the  previous argument cannot be used to infer that the deformation tensor vanishes identically. On the other hand, the same argument in Prop. 4.2 (iv) in \cite{PS} implies that the deformation tensor $\phi$ is bounded in  the $d d^c\t$-norm along each  leaf $L$. By holomorphicity of the component of $\phi|_L$ and Liouville Theorem,  this implies that  the components of $\phi|_L$ are constant along each leaf, hence identically vanishing if   $\f$  satisfies \eqref{biho} at the tangent spaces of $\pi^{-1}(\cS_o)$. \par
\smallskip
We have  now all ingredients to prove the three claims of the theorem. Let us start with (ii). In this case, for each model  $(\cM_o{(c)}, \t)$ the Remmert reduction $\pi: M_o{(c)} \to \cM_o{(c)}$ is the identity map and the assumption that   $\f$  satisfies \eqref{biho} at the  tangent spaces of $\pi^{-1}(\cS_o)$   coincides with the  condition that  $\f$ is a biholomorphism between the souls.  By the above discussion, this occurs if and only if $\f$ is a biholomorphism, proving that $(\cM_o{(c)}, \t)$ is soul rigid.\par
\smallskip
For the models considered in   (i) and  (iii), by the counterexamples in \S \ref{counterexamples} we know that  none of them is soul rigid. Nonetheless,  the above discussion shows that    when $M$ is of mixed type or  it has two ends and $c = + \infty$, if $\f$ is a biholomorphism at the blow ups of the souls,  then  it is a biholomorphism.  This  shows that in those cases, the model   $\cM_o{(c)}$ is soul semi-rigid, proving (iii) and the first claim of (i).   The  second claim in (i) is a consequence of the discussion  at the end of \S \ref{counterexamples}. 
\square
\par

\bigskip
\bigskip
\font\smallsmc = cmcsc8
\font\smalltt = cmtt8
\font\smallit = cmti8
\hbox{\parindent=0pt\parskip=0pt
\vbox{\baselineskip 9.5 pt \hsize=3.1truein
\obeylines
{\smallsmc
Morris Kalka
Mathematics Department
Tulane University
6823 St. Charles Ave.
New Orleans, LA 70118
USA
}\medskip
{\smallit E-mail}\/: {\smalltt kalka@math.tulane.edu}
}
\hskip -3truecm
\vbox{\baselineskip 8.5 pt \hsize=3.1truein
\obeylines
{\smallsmc
Giorgio Patrizio
Dip. Matematica e Informatica 
``U. Dini''
Universit\`a di Firenze
\& Istituto Nazionale di
Alta Matematica
``Francesco Severi''
}\medskip
{\smallit E-mail}\/: {\smalltt patrizio@math.unifi.it
}}
\hskip -2.5truecm
\font\smallsmc = cmcsc8
\font\smalltt = cmtt8
\font\smallit = cmti8
\hbox{\parindent=0pt\parskip=0pt
\vbox{\baselineskip 9.5 pt \hsize=3.1truein
\obeylines
{\smallsmc
Andrea Spiro
Scuola di Scienze e Tecnologie
Universit\`a di Camerino
Via Madonna delle Carceri
I-62032 Camerino (Macerata)
ITALY
}\medskip
{\smallit E-mail}\/: {\smalltt andrea.spiro@unicam.it}
}
}
}

%

\end{document}